\documentclass[a4paper,11pt,reqno]{amsart}
\usepackage{graphicx}
\usepackage{a4wide}
\usepackage{amssymb}
\usepackage{amsthm}
\usepackage{eucal}
\usepackage[pdftex,colorlinks]{hyperref}
\usepackage{pgf,pgfarrows,pgfnodes,pgfautomata,pgfheaps,pgfshade}
\usepackage{amsmath,amssymb,graphicx}
\usepackage{epsfig}

\numberwithin{equation}{section}

\newcommand{\R}{\mathbb R}

\newcommand{\Z}{\mathbb Z}

\newcommand{\T}{\mathbb T}

\newcommand{\ii}{\int\!\!\!\int } 

\newcommand{\be}{\begin{equation}}
\newcommand{\ee}{\end{equation}}
\newcommand{\ba}{\begin{eqnarray}}
\newcommand{\ea}{\end{eqnarray}}

\newtheorem{theorem}{Theorem}[section]
\newtheorem{proposition}[theorem]{Proposition}
\newtheorem{remark}[theorem]{Remark}
\newtheorem{lemma}[theorem]{Lemma}

\begin{document}

\title[Null controllability of a system of viscoelasticity]{Null controllability of a system of viscoelasticity\\
 with a moving control} 

\author{Felipe W. Chaves-Silva}
\address{BCAM � Basque Center for Applied Mathematics
Mazarredo 14, 48009 Bilbao, Basque Country, Spain}
\email{chaves@bcamath.org}

\author{Lionel Rosier}
\address{Institut Elie Cartan, UMR 7502 UdL/CNRS/INRIA,
B.P. 70239, 54506 Vand\oe uvre-l\`es-Nancy Cedex, France}
\email{Lionel.Rosier@univ-lorraine.fr}

\author{Enrique Zuazua}
\address{BCAM � Basque Center for Applied Mathematics
Mazarredo 14, 48009 Bilbao, Basque Country, Spain}
\address{Ikerbasque - Basque Foundation for
Science, Alameda Urquijo 36-5, Plaza Bizkaia,
48011, Bilbao, Basque Country, Spain}

\email{zuazua@bcamath.org}

\keywords{Wave equation with viscous and frictional damping; viscoelasticity; decoupling; moving control; null controllability; Carleman estimate}

\subjclass{}

\begin{abstract} 
In this paper, we consider the wave equation with both a viscous Kelvin-Voigt  and frictional  damping  
as a model of viscoelasticity in which we incorporate an internal control
with a moving support. We prove  the null controllability when the control region, driven by the flow of an ODE,  covers all the domain.
The proof is based upon the interpretation of the system as, roughly,  the coupling of a  heat equation with an ordinary differential equation (ODE). The presence of the ODE for which there is no propagation along the space variable makes the controllability of the system impossible when the control is confined into a subset in space that does not move. 
The null controllability  of the system with a moving control is established in using the observability of the adjoint system and some Carleman estimates for a coupled system of a parabolic equation and an ODE with the same singular weight, adapted to the geometry of the moving support of the control. This extends to the multi-dimensional case the results
 by P. Martin {\it et al.} on the one-dimensional case, employing $1-d$ Fourier analysis techniques.
\end{abstract}

\maketitle
\section{Introduction}
We are concerned with the controllability of the following model of viscoelasticity consisting of a wave equation with both viscous Kelvin-Voigt and frictional damping:
\ba
y_{tt}-\Delta y - \Delta y_t +b(x)y_t &=& 1_{\omega(t)}h,\qquad x\in \Omega, \ t\in (0,T), \label{A1}\\
y&=&0,\quad \qquad  \quad  x\in \partial \Omega,\ t\in (0,T), \label{A1bis}\\
y(x, 0)=y_0(x), \ y_t(x, 0) &= & y_1(x),\quad  \quad x\in \Omega. \label{A1bis2}
\ea
Here $\Omega$ is a smooth, bounded open set in $\R ^N$, $b\in L^\infty (\Omega)$ is a given 
function determining the frictional damping and $h=h(x,t)$ denotes the control.
To simplify the presentation and notation, and without loss of generality, the viscous constant has been taken to be the unit one $\nu = 1$. The same system could be considered with an arbitrary viscosity constant $\nu>0$ leading to the more general system
\begin{equation}\label{visco}
y_{tt}-\Delta y - \nu \Delta y_t +b(x)y_t = 1_{\omega(t)}h,
\end{equation}
but the analysis would be the same.

The control $h$ acting on the right hand side term as an external force is, for all $0<t<T$,  localized in a subset of $\Omega$. This fact is modeled by the multiplicative factor $1_{\omega(t)}$ which stands for the characteristic function of the set $\omega(t)$ that, for any $0<t<T$, constitutes  the support of the control, localized in a moving subset $\omega(t)$ of $\Omega$. 

Typically we shall consider control sets $\omega(t)$ determined by the evolution of a given reference subset $\omega$ of $\Omega$ through a smooth flow $X(x,t,0)$.

We consider the problem of null controllability. In other words, given a final time $T$ and initial data for the system $(y_0, y_1)$ in a suitable functional setting, we analyze the existence of a control $h=h(x, t)$ such that the corresponding solution satisfies the rest condition at the final time $t=T$:
$$
y(x, T) \equiv y_t(x, T) \equiv 0, \quad \hbox{ in } \Omega.
$$

One of the distinguished features of the system under consideration is that, for this null controllability condition to be fulfilled, the control needs to move in time. Indeed, if $\omega(t) \equiv \omega$  for all $0<t<T$, i.e. if the support of the control does not move in time as it is often considered,   the system under consideration is not controllable. This can be easily seen at the level of the dual observability problem. In fact, the structure of the underlying PDE operator  and, in particular, the existence of time-like characteristic hyperplanes, makes impossible the propagation of information in the space-like directions, thus making the observability inequality also impossible. This was already observed in the work by P. Martin {\it et al.} in \cite{MRR} in the $1-d$ setting.  There, for the $1-d$ model, it was shown that this obstruction could be removed by making the control move so that its support covers the whole domain where the equation evolves. 

More precisely, in \cite{MRR}, the $1-d$ version of the problem above was considered in the torus, with periodic boundary conditions,   $b\equiv 0$ and $\omega (t)=\{ x-t;\  x \in \omega \}$, i.e. 
\be
\label{A5}
y_{tt} - y_{xx} -y_{xxt}= 1_{\omega(t)}h(x,t), \quad x\in \T  .
\ee
Recall that  this system with boundary control, i.e. $h\equiv 0$ and the 
boundary conditions
\[
y(0,t)=0,\quad y(1,t)=g(t),
\]
$g=g(t)$ being the boundary control, 
fails to be spectrally controllable, because of  the existence of a limit point in the spectrum of the adjoint system \cite{RR}.
In the moving frame $x'=x+t$, \eqref{A5} 
may be written  as
\be
\label{A6}
z_{tt}  -2 z_{xt} -z_{txx} +z_{xxx}=a(x)h(x+t,t)
\ee 
where $z(x,t)=y(x+ t,t)$. In \cite{MRR} the spectrum of the adjoint system to \eqref{A6} was shown to be split into a hyperbolic part and 
a parabolic one. As a consequence, equation \eqref{A6} was proved to be null controllable in large time.   A similar result  was proved in \cite{RZ2012} for the 
Benjamin-Bona-Mahony equation
\[
y_t-y_{txx}+y_x+yy_x = a(x-ct)h(x,t),\qquad x\in \T.
\]
Once again this system turns out to be  globally controllable and exponentially stabilizable in $H^1(\T)$ for any $c\ne 0$.
But, as noticed in \cite{Micu}, the linearized equation fails to be spectrally controllable with a control supported in a fixed domain. 

As mentioned above, in both cases, the lack of controllability of these systems with immobile controls is due to the fact that the underlying PDE operators exhibit the presence of time-like characteristic lines thus making propagation in the space-like directions impossible. By the contrary, when analyzing the problem in a moving frame, the characteristic lines are oblique ones in $(x, t)$, thus facilitating propagation properties.

The main goal of this paper is to extend the $1-d$ analysis in \cite{MRR} to the multi-dimensional case. This can not be done with the techniques in \cite{MRR} based on Fourier analysis.  Our approach is rather inspired on the fact that system (\ref{A1})-(\ref{A1bis}) can be rewritten as a system  coupling a parabolic equation with an ordinary differential equation (ODE). The presence of this ODE, in the case of a fixed support of the control, independent of $t$,  is responsible for the lack of controllability of the system, due to the absence of propagation in the space-like direction.  Letting the control move introduces an effect similar to adding a transport term in the ODE but keeping the control immobile,  thus changing the structure of the system into a parabolic-transport coupled one. This new system turns out to be controllable under the condition that all characteristics of the transport equation enter within the control set in the given control time, a condition that is reminiscent of the so-called Geometric Control Condition in the context of the wave equation (see \cite{BLR}).

The approach in \cite{MRR} would suggest to do the following splitting of \eqref{A1}: 
\ba
v_t -\Delta v &=& 1_{\omega(t)} h + (1-b)(v-y) \label{A11intro}\\
y_t+y &=& v.\label{A12intro}
\ea

However,  the splitting can be performed in an alternative manner as follows:
\ba
y_t -\Delta y + (b-1) y  &=& z \label{A13bis}\\
z_t+z &=& 1_{\omega(t)} h +  (b-1)y.\label{A14bis}
\ea
It can  easily  be seen that $y$ solves equation (\ref{A1}) if, and only if, it is the first component of the solution of system (\ref{A13bis})-(\ref{A14bis}). 


Our analysis of the Carleman inequalities for these systems is  analog to that in  \cite{AT}
for a system of thermoelasticity coupling the heat and the wave equation. The key in \cite{AT} and in our own analysis is to use the same weight function both for the Carleman inequality of the heat and the hyperbolic model. In \cite{AT}, since dealing with the wave equation, rather strong geometric conditions were needed on the subset where the control or the observation mechanism acts. In our case, since we are considering the simpler transport equation, the geometric assumptions will be milder, consisting mainly on a characteristic condition ensuring that all characteristic lines of the transport equation intersect the control/observation set. This suffices for the Carleman inequality to hold for the transport equation and is also sufficient for the heat equation that it is well-known to be controllable/observable from any open non-empty subset of the space-time cylinder where the equation is formulated.

It is important to mention that, as far as we know, all the Carleman inequalities for the heat equation  available in the literature are done for the case where the control region is fixed. In the case we are dealing, the control region  is moving in time. Therefore, a new Carleman inequality must be proved in this framework. The proof of a Carleman inequality for the heat equation when the control region is moving  is one of the novelties of this paper.

In order to state  the main result of this paper we first  describe precisely the class of moving trajectories for 
the control for which our null controllability result will hold.

\null

\textbf{Admissible trajectories}:  In practice  the trajectory of the control  can be taken to be determined by the flow $X(x,t,t_0)$ generated by some vector 
 field $f\in C([0,T];W^{2,\infty}(\R ^N;\R ^N ))$, i.e. $X$ solves
\be
\left\{ 
\begin{array}{l}
\displaystyle \frac{\partial X}{\partial t} (x , t , t_0)= f ( X ( x,t,t_0),t),\\[3mm]
X(x, t_0 , t_0)=x. 
\end{array}
\right.
\label{Xflow}
\ee
For instance, any translation of the form:
\be
\label{translation}
X(x,t,t_0)=x+ \gamma (t) -\gamma (t_0),
\ee
where $\gamma \in C^1([0,T];\R ^N)$, is admissible. (Pick $f(x,t)=  \dot \gamma (t)$).

We assume that there exist a bounded, smooth,  open set $\omega _0\subset \R ^N$, a curve $\Gamma \in C^\infty ([0,T];\R ^N)$, 
and two times $t_1,t_2$ with $0\le t_1 < t_2 \le T $ such that: 
\ba
&& \Gamma (t) \in X(\omega _0,t,0)\cap \Omega , \quad \forall t\in [0,T]; \label{A3a} \\[3mm]
&&\overline{\Omega} \subset \displaystyle\cup_{t\in [0,T]} X(\omega _0,t,0) =\{ X(x,t,0);\ \ x\in \omega _0,\ t\in [0,T]\} ; \label{A3b} \\[3mm]
&&\Omega \setminus \overline{X(\omega _0,t,0) } \text{ is nonempty and connected for } t\in [0,t_1]\cup [t_2,T]; \label{A3c}\\ [3mm]
&&\Omega \setminus \overline{X(\omega _0,t,0) } \text{ has two (nonempty) connected components  for } t\in (t_1,t_2); \label{A3d}\\[3mm]
&&\forall \gamma\in C([0,T]; \Omega ) , \ \exists t\in [0,T],\quad \gamma (t) \in X(\omega _0,t,0). \label{A3e}
\ea


  
The main result in this paper is as follows.
\begin{theorem}
\label{thm1}
Let $T>0$, $X(x,t,t_0)$ and  $\omega_0$ be as in \eqref{A3a}-\eqref{A3e}, and let $ \omega$ be  any open set in $\Omega$ such that $ \overline{\omega}_0\subset \omega$. 
Then for all $(y_0,y_1)\in L^2 (\Omega)^2$  with $y_1 - \Delta y_0\in L^2(\Omega )$, there exists a
function $h\in L^2(0,T;L^2(\Omega ))$ for which the solution of 
\ba
&&y_{tt}-\Delta y - \Delta y_t +b(x)y_t= {\bf 1}_{X( \omega  ,t,0) } (x)   h, \qquad (x,t) \in \Omega \times (0,T),\label{A31}\\
&&y(x,t)=0, \qquad (x,t)\in \partial \Omega \times (0,T),\label{A32}\\ 
&&y(.,0)=y_0,\ y_t(.,0)=y_1,\label{A33}
\ea
fulfills $y(.,T)=y_t(.,T)=0$.
\end{theorem} 

A few remarks are in order in what concerns the functional setting of this model:
\begin{itemize}
\item Viewing the system of viscoelasticity under consideration as a damped wave equation, a natural functional setting would be the following: For data in $H^1_0(\Omega) \times L^2(\Omega)$ and, say, right hand side term of (\ref{A1}) in $L^2(0, T; H^{-1}(\Omega))$, there exists an unique solution $y \in C([0, T]; H^1_0(\Omega) \cap C^1([0, T]; L^2(\Omega))$. Furthermore $y_t \in L^2(0, T; H^1_0(\Omega))$. The latter is an added integrability/regularity property of the solution that is due to the strong damping effect of the system. This can be seen naturally by considering the energy of the system
$$
E(t) = \frac{1}{2} \int_\Omega [�|y_t|^2 + |\nabla y |^2 ] dx,
$$
that fulfills the energy dissipation law
$$
\frac{d}{dt} E(t) = -\int_\Omega [�|\nabla y_t|^2 + b(x) |y_t|^2 ] dx + \int_{\omega(t)} h y_t dx.
$$

\item We can also solve (\ref{A13bis})-(\ref{A14bis}) so that $y$, solution of the heat equation, lies in the space $y\in C([0, T]; H^1_0(\Omega)) \cap L^2(0, T; H^2(\Omega))$ and $z$, solution of the ODE, in $C([0, T]; L^2(\Omega) )$. This can be done provided $(y_0, y_1-\Delta y_0)\in H^1_0(\Omega) \times L^2(\Omega)$. The functional setting is not exactly the same but this is due to the fact that, in some sense, in one case we see the system as a perturbation of the wave equation, while, in the other one,  as a variant of the heat equation.

\item From a control theoretical point of view it is much more efficient to analyze the system in the second setting, as a perturbation of the heat equation, through the coupling with the ODE or, after changing variables, with a transport equation. If we view the system of viscoelasticity as a perturbation of the wave equation, standard hyperbolic techniques such as multiplier, Carleman inequalities or microlocal tools do not apply since, actually, the viscoelastic term determines the principal part of the underlying PDE operator and cannot be viewed as a perturbation of the wave dynamics.

\end{itemize}

The analysis is particularly simple in the special case where   $b\equiv 1$. Indeed, in that case, the system \eqref{A11intro}-\eqref{A12intro} 
(with a second control incorporated in the ODE) takes the following cascade form
\begin{eqnarray}
v_t -\Delta v &=& 1_{\omega ( t ) } \tilde h, \label{A34}\\
y_t + y &=& 1_{\omega ( t ) }  \tilde k +  v,   \label{A35}
\end{eqnarray}
where the parabolic equation \eqref{A34} is {\em uncoupled}. This system will be investigated  in a separate section (Section \ref{uncoupled}) since some of the basic ideas allowing to handle the general case emerge already in its analysis.  Note that, in this particular case, roughly, one can first control the heat equation by a suitable control $\tilde h$ and then, once this is done, viewing $v$ as a given source term,  control the transport equation by a convenient  $\tilde k$. This case is also important because the only assumption needed to prove Theorem \ref{thm1} in this case is 
\eqref{A3b}  (i.e. we don't assume that \eqref{A3a} and \eqref{A3c}-\eqref{A3e} are satisfied). 

 In this particular case $b\equiv 1$ a similar argument can be used with the second decomposition.

The paper is organized as follows. Section \ref{uncoupled} is devoted to address the particular case $b \equiv 1$.
In Section \ref{Examples} we give some examples showing  the importance of the taken  assumptions on the trajectories. 
In Section \ref{Nullcontrol}, we go back to the general system  \eqref{A13bis}-\eqref{A14bis}. We prove that this system is null controllable in $L^2(\Omega )^2$ 
(see Theorem \ref{thm3}) by deriving the 
observability of the adjoint system from two Carleman estimates with the same singular weight, adapted to the flow determining the moving control. The details of the construction of the  weight function based on \eqref{A3a}-\eqref{A3e} are given in Lemma \ref{weight}. Theorem \ref{thm1} is then a direct consequence of Theorem \ref{thm3}. 
We finish this paper with two further sections devoted to comment some closely related issues and open problems

\section{Analysis of the decoupled cascade system}\label{uncoupled}
In this section, we give a proof of Theorem \ref{thm1} in the special situation when $b\equiv 1$, and $\omega (t)= X (\omega  _0, t, 0)$, where 
$X$ is given by \eqref{Xflow} for some $f\in C(\R ^+ ;W^{2,\infty}(\R ^N;\R ^N ))$.

\null

As we said before, we will prove Theorem \ref{thm1} in the case $b\equiv 1$ by proving  a null controllability result for  the decomposition  \eqref{A34}-\eqref{A35}. The idea of the proof is as follows. We take some appropriate $ 0 < \epsilon < T$ and then  drive the solution of  the heat equation \eqref{A34} to zero in time $\epsilon$ by means of a control $\tilde{h}$. Next, we let equation \eqref{A35}  evolves freely in $[0,\epsilon]$, i.e., $\tilde{k}\equiv 0$, and then we control this equation by means of a smooth 
control $\tilde{k}$ in the time interval $[\epsilon, T]$. This gives the null controllability of the system \eqref{A34}-\eqref{A35} in the whole time interval $[0,T]$.

\null

\textbf{Proof of Theorem \ref{thm1} in the case $b\equiv 1$}.

\null 

Suppose \eqref{A3b}  is satisfied and let
\be
\label{D1}
T_0 =\inf\{ T>0;\  \overline{\Omega}  \subset \cup_{0\le t\le T }X( \omega _0,t,0 ) \} .
\ee

Pick any $T>T_0$, and pick some $\varepsilon \in (0,T-T_0)$ and some nonempty open set $\omega _{-1}\subset \omega_0$ such that 
\ba
\overline{\Omega } &\subset& \cup_{\varepsilon \le t\le  T} X(\omega _0, t, 0),\label{time}\\
\omega  _{-1}&\subset& X(\omega _0,t,0)\quad \forall t\in (0,\varepsilon). 
\label{omega1}
\ea
Let  $T'=T- \varepsilon$, and pick any $(v_0,y_0)\in L^2(\Omega )^2$. Then, it is well known (see \cite{Imanuvilov})  that there exists some control input 
$ h\in L^2(0,\varepsilon ; L^2(\Omega ))$ such that the solution $v=v(x,t)$ of 
\ba
v_t - \Delta v  = {\bf 1}_{\omega _{-1}} h ,&& x\in \Omega ,\  t\in (0, \varepsilon ) , \\
v(x,0)=v_0(x),&&x\in \Omega, 
\ea 
satisfies 
\[
v(x,\varepsilon ) =0, \qquad x\in \Omega .
\] 
Set 
\begin{eqnarray*}
\tilde h(x,t) = {\bf 1}_{\omega  _{-1}}(x) h(x,t),&&   x\in \Omega ,\  t\in ( 0,\varepsilon ), \\
\tilde h(x,t) = 0,&&x\in \Omega ,\  t\in ( \varepsilon ,T), \\
\tilde k(x,t)=0,&& x\in \Omega , \  t\in (0,\varepsilon )
\end{eqnarray*}
Then the solution $v$ of 
\ba
v_t - \Delta v  = {\bf 1}_{X(\omega _0 ,t,0)} (x)\tilde h&& x\in \Omega ,\  t\in (0, T) , \\
v(x,0)=v_0(x),&&x\in \Omega, 
\ea 
satisfies $v(x,t)=0$ for $t\in [\varepsilon , T]$. 
We claim that the system
\ba
&&y_t +y = {\bf 1}_{X(\omega  _0,t,0)}(x)  k(x,t) , \qquad x\in \Omega , \  t\in (\varepsilon ,T),\label{D7}\\
&&y(x,\varepsilon )=y_0(x),\label{D8}
\ea
is exactly controllable in $L^2(\Omega )$ on the time interval  $(\varepsilon , T)$. By duality, this is equivalent to proving that the corresponding observability inequality
\be
\label{D4}
\int_{\Omega } |q_0(x)|^2 dx \le C \int_{\varepsilon}^{T}\!\!\!\int_{\Omega } {\bf 1}_{X(\omega _0 ,t,0)} (x)\ |q(x,t)|^2 dxdt 
\ee 
is fulfilled with a uniform constant $C>0$ for all solution $q$ of the adjoint system 
\ba
&&-q_t + q =0, \qquad x\in \Omega ,\  t\in (\varepsilon , T),\label{D10}\\
&&q(x,T)=q_0(x). \label{D11}
\ea
Since the solution of \eqref{D10}-\eqref{D11} is given by $q(x,t)=e^{t-T}q_0(x)$, we have that 
\begin{equation}
\int_{\varepsilon} ^{T} \!\!\! \int_{\Omega } {\bf 1}_{X(\omega _0 ,t,0)} (x)\  |q(x,t)|^2 dxdt  \ge 
e^{2(\varepsilon -T)} \int_{\Omega}  |q_0(x)|^2 (\int_{\varepsilon}^{T}  {\bf 1}_{X(\omega _0 ,t,0)} (x)\ dt) dx .\label{D25}
\end{equation} 
From \eqref{D1} and the smoothness of $X$, we see that for all $x\in \overline{\Omega}$, there is some $t_0\in (\varepsilon ,T)$, and some 
$\delta >0$ such that for any $y\in B(x,\delta )$ and any $t\in (\varepsilon ,T)\cap (t_0-\delta, t_0+\delta )$ we have $y\in X(\omega _0, t,0)$. 
From the compactness of $\overline{\Omega}$, we see that there exists a number $\delta _0 >0$ such that
\[
\int_{\varepsilon} ^{T} {\bf 1}_{X(\omega _0,t,0)  } (x) dt >\delta _0, \qquad \forall x\in \overline{\Omega }  .  
\] 
Combined with \eqref{D25}, this yields \eqref{D4}. 
Thus, \eqref{D7}-\eqref{D8} is exactly controllable in $L^2(\Omega )$ on $(\varepsilon , T)$ with some controls
$k\in C([\varepsilon ,T];L^2(\Omega ))$. 
Let $y_1(x)=e^{-\varepsilon} y_0(x) + \int_0^\varepsilon e^{s-\varepsilon } v(x,s)ds $. Extend $\tilde k$ to $(0,T)$ so that 
$\tilde k\in L^2(0,T ; L^2(\Omega))$ and the solution of 
\ba
&&y_t +y = {\bf 1}_{X(\omega  _0,t,0)}(x)  \tilde k , \qquad x\in \Omega , \  t\in (\varepsilon , T),\\
&&y(x,\varepsilon)=y_1(x),
\ea
satisfies $y(.,T)=0$.  Thus the control $(\tilde h,\tilde k)$ steers the solution of \eqref{A34}-\eqref{A35} from $(v_0,y_0)$ at $t=0$ to  $(0,0)$ at $t=T$. 
Applying the operator $\partial _t-\Delta$ in each side of \eqref{A35} results in 
\be 
y_{tt}-\Delta y - \Delta y_t +y_t= {\bf 1}_{X( \omega _0  ,t,0) } (x)   \tilde h +(\partial _t -\Delta ) [ {\bf 1}_{X( \omega  _0 ,t,0) } (x)\tilde k ],
\label{final}
\ee
This proves Theorem \ref{thm1}, except the fact that  the control does not live in $L^2(0,T;L^2(\Omega ))$. Assume now
that $(v_0,y_0)\in L^2(\Omega ) \times [H^2(\Omega )\cap H^1_0(\Omega )] $. To get a control
$\tilde k\in L^2(0,T;H^2(\Omega ))$, it is sufficient to replace ${\bf 1} _{\omega  (t)}$ by $a(X(x, 0,t))$ in \eqref{A35}, where   $a$ is a function satisfying
\begin{eqnarray*}
&& a\in C_0^\infty (\omega ),\\
&& a(x)=1\quad \forall x\in \omega_0.
\end{eqnarray*}
The proof is completed by  showing  the observability inequality
\[
||q_0||^2_{X'} \le C \int_{\varepsilon} ^{T}  ||(t-\varepsilon) a(X(., 0,t)) q(.,t)||^2_{X'}dt
\] 
for the solution $q$ of system \eqref{D10}-\eqref{D11}, 
where $X=H^2(\Omega )\cap H^1_0(\Omega )$ and $X'$ stands for its dual space. 
This can be done as in \cite[Proposition 2.1]{RZ2010}. Next, using the HUM operator, we notice that $\tilde k  \in C^1([\varepsilon ,T];X)$ with $\tilde k(.,\varepsilon)=0$, 
since $q\in C^1( [ \varepsilon , T];X' )$ for any $q_0\in X'$. Thus, with this small change,  the right hand side term in \eqref{final}  can be written
${\bf 1} _{X(\omega , t,0)}u(x,t)$, where $u\in L^2(0,T;L^2(\Omega ))$.

\begin{remark}
Observe that the situation when $X(\omega  _0,t,0)$ moves as in 
Figure \ref{fig5}, Figure \ref{fig3}  or in Figure \ref{fig2} (see below) is admissible in the case when $b \equiv 1$. 
\end{remark}

\section{Examples}\label{Examples}
In this section, we provide some geometric examples to illustrate the assumptions \eqref{A3a}-\eqref{A3e}. 
We use simple shapes (like rectangles) just for convenience. 

\begin{itemize}
\item Figure \ref{fig1} shows how a control region should move in order to satisfy conditions \eqref{A3a}-\eqref{A3e}. 
\begin{figure*}[!t]
\centerline{\includegraphics[width=150mm]{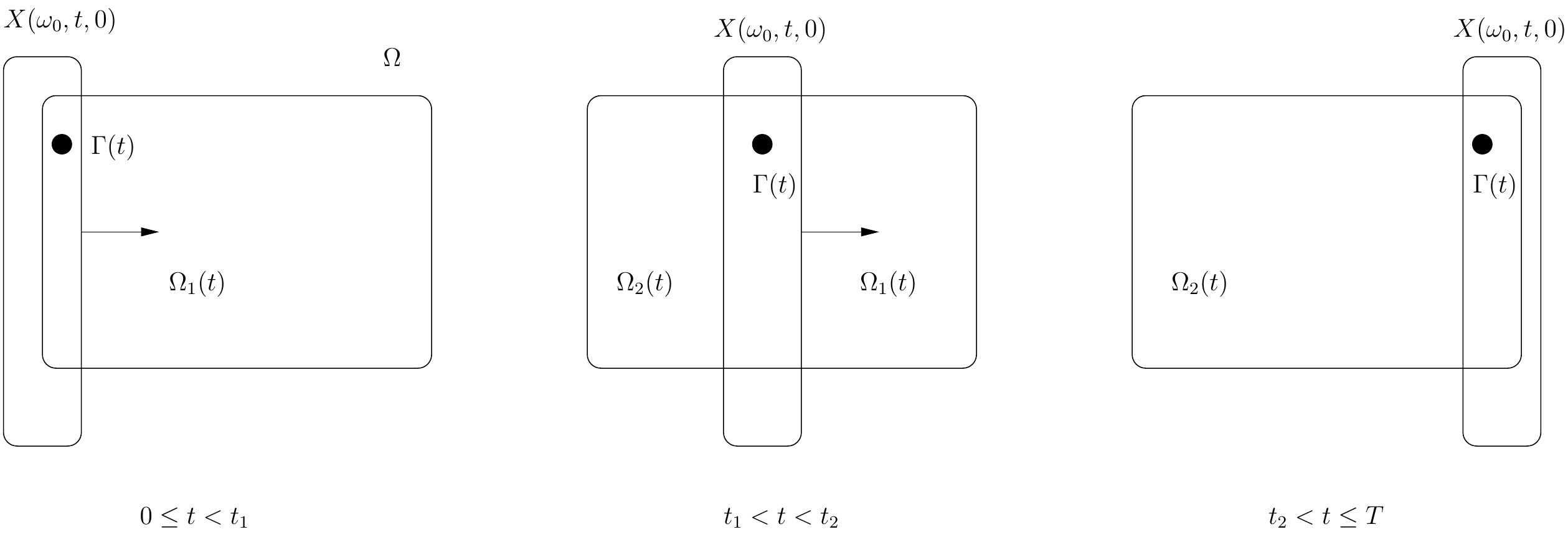}}
\caption{Example for which conditions \eqref{A3a}-\eqref{A3e} are satisfied.}
\label{fig1}
\end{figure*}
%
\item 
Figure \ref{fig5} depicts a situation for which Theorem \ref{thm1} cannot  be applied, except in the case when  $b \equiv 1$, as condition \eqref{A3d} fails. 
\begin{figure*}[!t]
\centerline{\includegraphics[width=150mm]{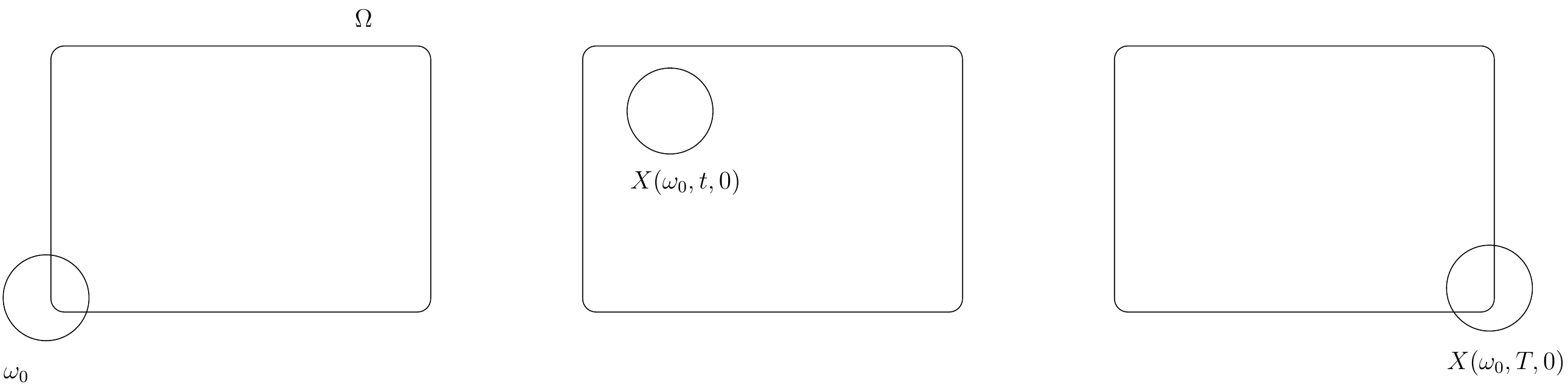}}
\caption{Example for which condition \eqref{A3d} fails.}
\label{fig5}
\end{figure*}
%
\item
\begin{figure*}[!t]
\centerline{\includegraphics[width=150mm]{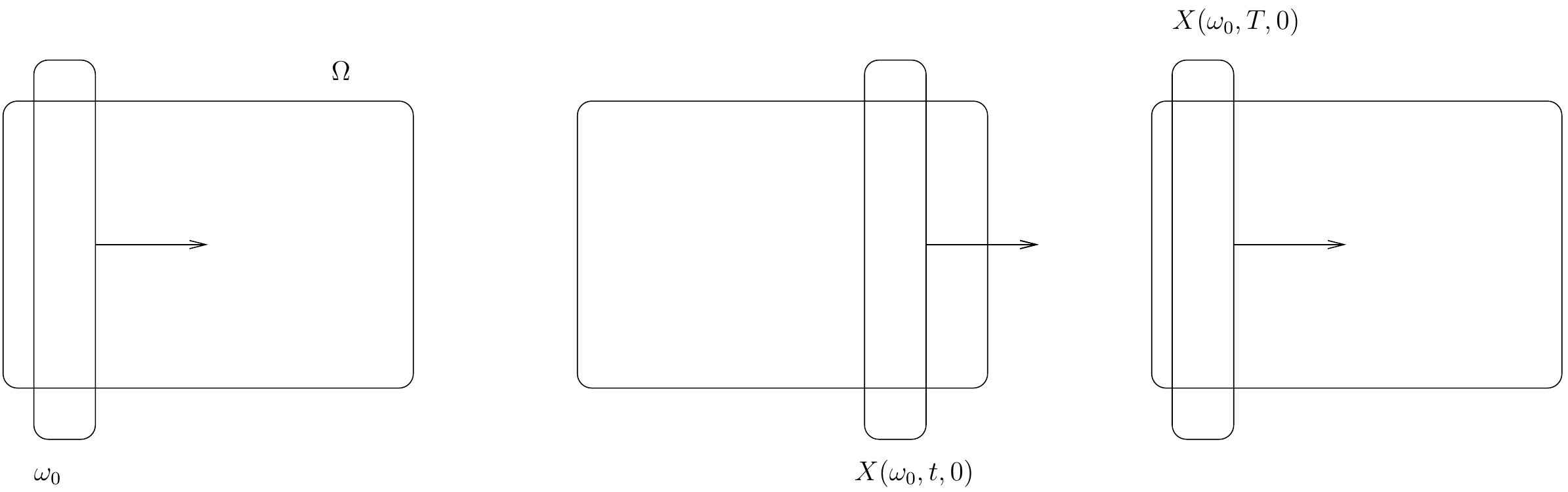}}
\caption{Example for which condition \eqref{A3c} fails.}
\label{fig6}
\end{figure*}
In Figure \ref{fig6}, we modify the example given in Figure \ref{fig1} by shifting the time. Theorem \ref{thm1} cannot be applied as it is, since
\eqref{A3c} fails. However,  the conclusion of Theorem \ref{thm1} remains valid. Indeed, assume that $\Omega \setminus \overline{\omega (t) }$ has 
two connected components (resp. one) for $t\in [0,\tau _1)\cup (\tau _2,T]$ (resp. for $t\in [\tau _1,\tau _2]$).  Assume that the ``jump'' of $\omega (t)$  occurs at $t=\tau _3$, with 
$\tau _1<\tau _3<\tau _2$. 
Let 
\ba
&&{\mathcal O}_1  := \cup _{0\le t\le \tau _3}\ \omega  (t)  ,\\ 
&&{\mathcal O}_2 := \cup _{\tau _3 \le t \le T}\ \omega  (t)   
\ea
and let $\eta  \in C^\infty (\Omega ;[0,1])$ be such that 
\ba
&&\textrm{supp} (\eta )\subset {\mathcal O}_1,\\
&&\textrm{supp} (1-\eta )\subset {\mathcal O}_2,\\
&&\textrm{supp} (\nabla \eta ) \subset \omega _0.
\ea
Then, applying the Carleman estimate in Lemma \ref{lem3} to $(p_1,q_1)=\eta \big( X(x,0,t)\big) (p,q)$  in $\Omega \cap \{ \eta >0 \} $
on the time interval $[0,\tau _3]$, and to $(p_2,q_2)=\big( 1-\eta ( X(x,0,t) )  \big) ( p,q)$  in $\Omega \cap \{ \eta  < 1 \} $
   on the time interval $[\tau _3,T]$, we can easily prove the observability inequality \eqref{observability}.  
\item Figure \ref{fig3} shows that the assumption \eqref{A3a}, which is needed to construct the weight function $\psi$ in  Lemma \ref{weight} 
 cannot be replaced by the simpler condition
\[
X(\omega _0,t,0)\cap \Omega \ne \emptyset , \quad \forall t\in[0,T].
\]
\begin{figure*}[!t]
\centerline{\includegraphics[width=150mm]{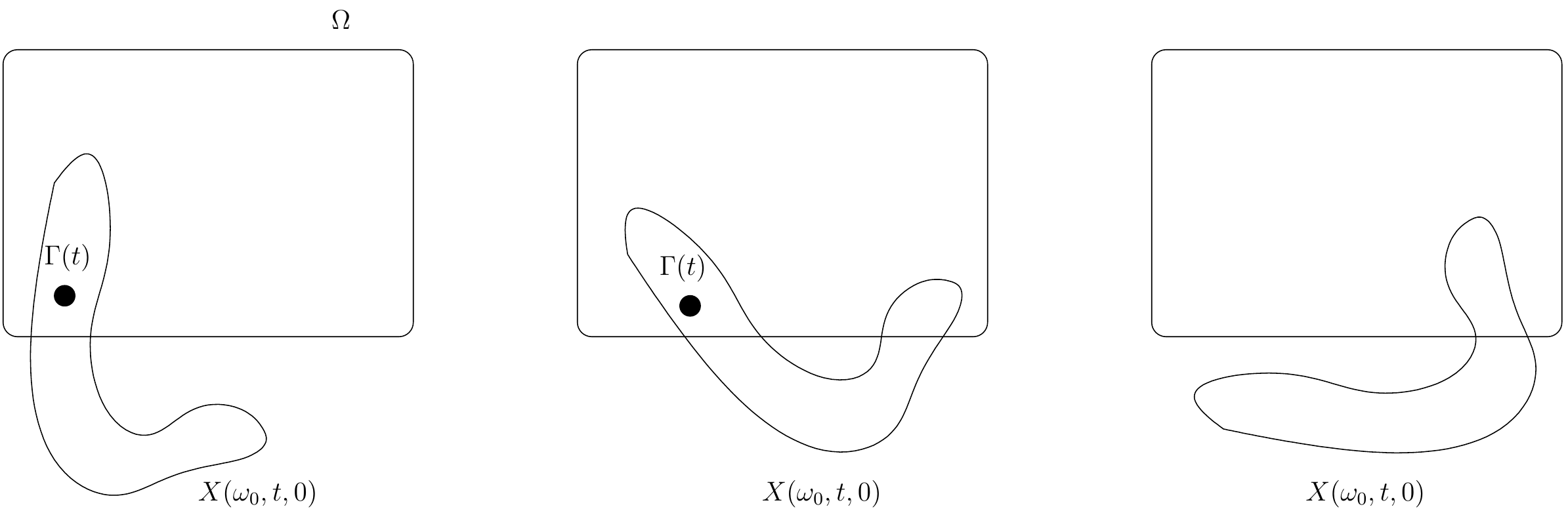}}
\caption{Example showing that $X(\omega _0,t,0)\cap \Omega \ne \emptyset$ $\forall t\in[0,T]$ does not imply \eqref{A3a}.}
\label{fig3}
\end{figure*}
%
\item  Figure \ref{fig2} shows that the assumption \eqref{A3e}, which is also needed to construct the weight function $\psi$ in Lemma \ref{weight}, does not result 
from the other assumptions \eqref{A3a}-\eqref{A3d}.

\begin{figure*}[!t]
\centerline{\includegraphics[width=150mm]{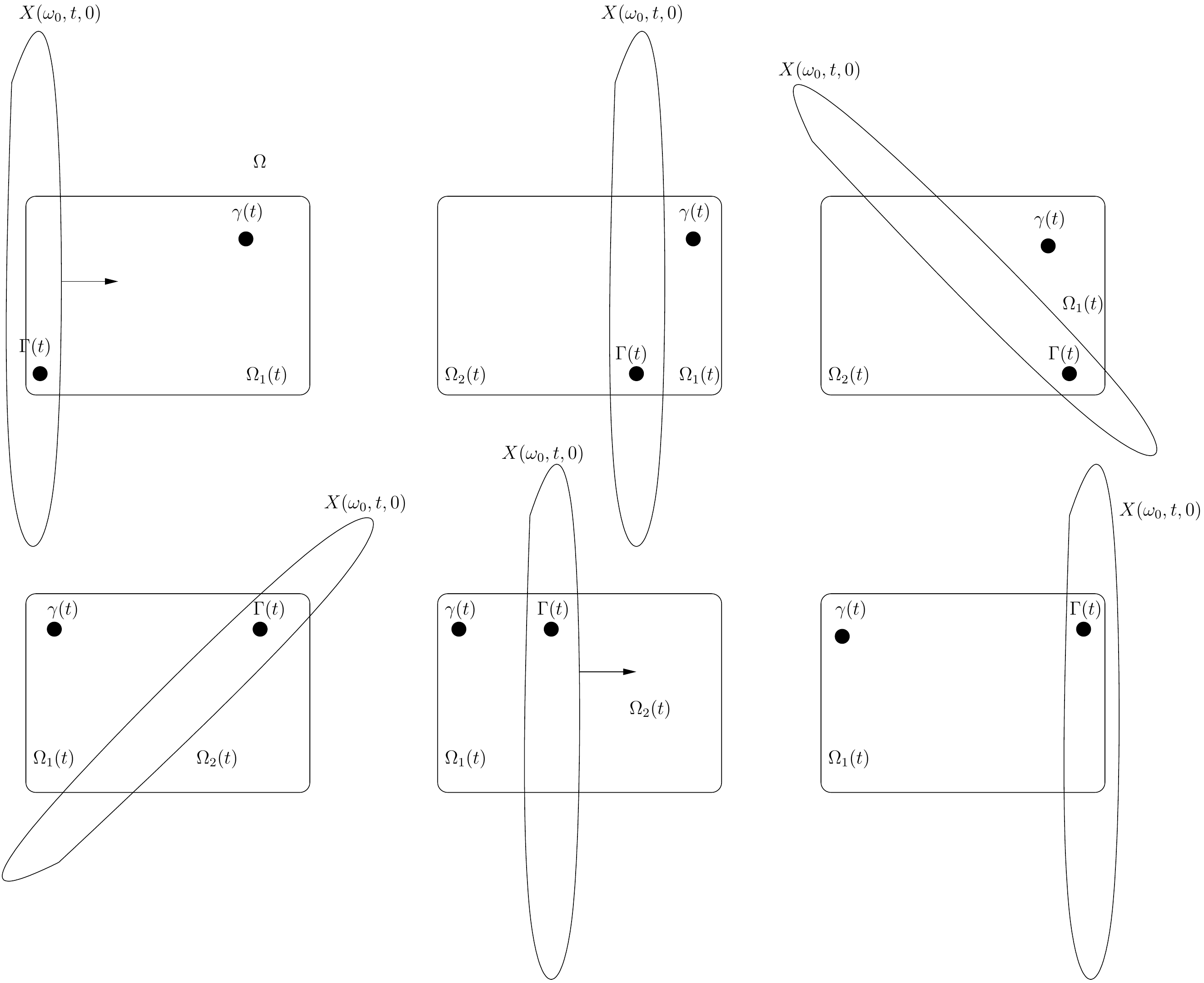}}
\caption{Example showing that \eqref{A3a}-\eqref{A3d} does not imply \eqref{A3e}.}
\label{fig2}
\end{figure*}
%
\end{itemize}

\section{Null controllability of system \eqref{A34}-\eqref{A35}.}\label{Nullcontrol}
 
In this section we proof Theorem \ref{thm1}.  Using  decomposition \eqref{A13bis}-\eqref{A14bis}, it is easy to see that the null controllability 
of \eqref{A1}-\eqref{A1bis2}  turns out to be equivalent to the null controllability  of  the system

\begin{eqnarray}
 y_t -\Delta y  + (b(x)-1)y &=& z,\qquad (x,t)\in \Omega \times (0,T), \label{E1-1}\\
z_t+z &=& {\bf 1} _{X (\omega ,t,0)} (x) h+   (b(x)-1)y,\qquad (x,t)\in \Omega \times (0,T),   \label{E2-1}\\
y(x,t)&=&0,\qquad (x,t)\in \partial \Omega \times (0,T),\label{E3-1} \\
z(x,0)&=&z_0(x), \qquad x\in\Omega ,\label{E4-1}\\
y(x,0)&=& y_0(x), \qquad x\in \Omega. \label{E5-1}
\end{eqnarray}

More precisely, Theorem \ref{thm1} is a direct consequence of the following result.
\begin{theorem}
\label{thm3}
Let $T$, $X(x,t,t_0)$ and  $\omega_0$ be as in \eqref{A3a}-\eqref{A3e}, and let $\omega$ be as in Theorem \ref{thm1}. 
Then for all $(y_0,z_0)\in L^2(\Omega)^2$, there exists a control
function $h\in L^2(0,T;L^2(\Omega  ))$ for which the solution $(y,z)$ of \eqref{E1-1}-\eqref{E5-1}
satisfies $y(.,T)=z(.,T)=0$.
\end{theorem}


From now on we concentrate in the proof of  Theorem \ref{thm3}.

 It is well-known (see \cite{FI} ) that  Theorem \ref{thm3} is equivalent  to prove an observability inequality for the  adjoint system of \eqref{E1-1}-\eqref{E5-1}, namely
\ba
-p_t -\Delta p +(b(x) -1 )p &=& (b(x) -1 ) q, \qquad (x,t)\in \Omega \times (0,T), \label{B1000-1}\\
-q_t +q &=& p, \qquad (x,t)\in \Omega \times (0,T), \label{B2000-1}\\
p(x,t) &=& 0, \qquad (x,t) \in \partial \Omega \times (0,T), \label{B3000-1}\\
p(x,T)&=&p_0(x), \qquad x\in\Omega ,\label{B4000-1}\\
q(x,T)&=& q_0(x), \qquad x\in \Omega. \label{B5000-1}
\ea

In fact,  one can show that Theorem \ref{thm3} is equivalent to the following:
\begin{proposition}
\label{prop1} 
Let $T$, $X$, $\omega _0$ and $\omega $ be as in Theorem \ref{thm1}. Then there exists a constant $C>0$ such that for all $(p_0,q_0)\in L^2(\Omega  )^2$, the solution 
$(p,q)$ of \eqref{B1000-1}-\eqref{B5000-1} satisfies
\be
\label{observability}
\int_{\Omega} [|p(x,0)|^2 + |q(x,0)|^2] dx \le C \int_0^T\!\!\!\int_{ X(\omega ,t,0)}   |q(x,t)|^2\ dxdt.
\ee
\end{proposition}


\noindent
{\em Proof of Proposition \ref{prop1}.}
Inspired in part by \cite{AT} (which was concerned with a heat-wave system\footnote{See also \cite{EGGP} for some Carleman estimates for a coupled system of parabolic-hyperbolic equations.}), 
we shall establish some Carleman estimates for the (backward) parabolic equation \eqref{B1000-1} and the ODE \eqref{B2000-1} with the {\em same singular weight}.

\null

For a better comprehension, the proof  will be divided into two steps as follows:

\null

\textit{Step 1.} We apply  suitable  Carleman estimates for the parabolic equation \eqref{B1000-1} and the ODE \eqref{B2000-1}, 
with the same weights and with a moving control region.

\null
 
\textit{Step 2.}  We estimate a local integral of  $p$ in terms of a local integral of $q$ and some small order terms. 
Finally, we  combine all the estimates obtained in the first step and derive the desired Carleman inequality.

\null

The basic weight function we need in order  to prove such inequalities is given by the following Lemma.
\begin{lemma}
\label{weight}
Let $X$, $\omega_0$ and $\omega $ be as in Theorem \ref{thm1}, and let $\omega_1$ be  a nonempty open set in $\mathbb{R}^N$ such that
\be
 \overline{\omega _0}\subset \omega _1, \quad \overline{\omega _1} \subset \omega.  \label{A3-x}
\ee

Then there exist a number $\delta \in (0,T/2)$ and 
 a function $\psi \in C^\infty ( \overline{\Omega} \times [0,T])$ such that
\ba
\nabla \psi(x,t) \ne 0,&&\quad t\in[0,T],\ x\in \overline{\Omega} \setminus X(\omega _1,t,0),  \label{P1}\\
\psi _t (x,t) \ne 0,  &&  \quad  t\in[0,T],\ x\in \overline{\Omega}\setminus X(\omega _1,t,0),  \label{P2}\\
\psi _t (x,t)  >0,  &&  \quad  t\in[0,\delta ],\ x\in \overline{\Omega}\setminus X(\omega _1,t,0), \label{P3} \\
\psi _t  (x,t) <0,  &&  \quad  t\in[T-\delta ,T ],\ x\in \overline{\Omega}\setminus  X(\omega _1,t,0) ,\label{P4} \\
\frac{\partial \psi}{\partial n}(x,t) \le 0,  && \quad  t\in [0,T ],\ x\in \partial \Omega , \label{P5}\\
\psi (x,t) >\frac{3}{4}||\psi ||_{L^\infty (\Omega \times (0,T)) }, &&\quad  t\in [0,T ],\ x\in \overline{ \Omega} .\label{P6}
\ea
\end{lemma}

The proof of Lemma \ref{weight} will be given in Appendix \ref{aa}.
\begin{figure*}[!t]
\centerline{\includegraphics[width=100mm]{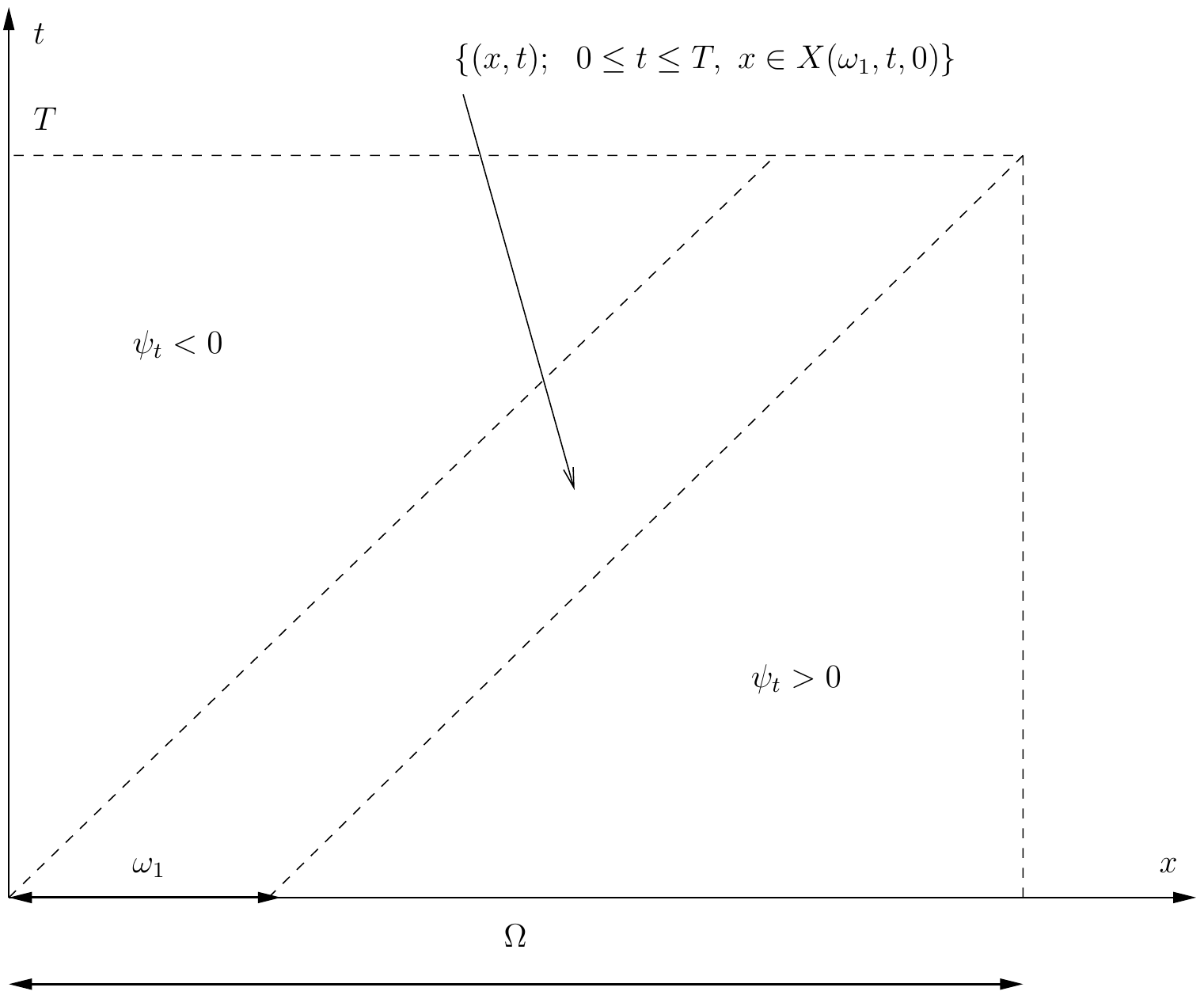}}
\caption{Sign of the time derivative of $\psi$.}
\label{fig4}
\end{figure*}


Next, we pick a function  $g\in C^\infty (0,T)$ such that 
\[
g(t) = \left\{ 
\begin{array}{ll}
\frac{1}{t} \qquad &\text{for } 0<t<\delta /2 , \\
\text{\rm strictly decreasing}\qquad &\text{for }  0<t \le  \delta ,\\
1 \qquad &\text{for } \delta \le t \le \frac{T}{2},\\
g(T-t)\qquad &\text{for }   \frac{T}{2} \le t < T
\end{array}
\right.
\]
 and define the weights
\[
\begin{array}{rll}
\varphi (x,t) &= g(t) (e^{\frac{3}{2} \lambda ||\psi ||_{L^\infty}} - e^{\lambda \psi (x,t)}), \qquad &(x,t)\in \Omega \times (0,T),\\
\theta (x,t) &= g(t) e^{ \lambda \psi (x,t)},\qquad &(x,t)\in \Omega  \times (0,T),
\end{array}
\] 
where $||\psi ||_{L^\infty}= ||\psi ||_{L^\infty (\Omega \times (0,T)) }$ and $\lambda >0$ is a parameter. 

\null


\textit{Step 1.} {\sc Carleman estimates with the same weight.}

In this step we apply a Carleman inequality for the heat-like equation \eqref{B1000-1} and a   Carleman inequality for the ODE \eqref{B2000-1}, both with the same weight. 
We combine such inequalities and  obtain a global estimation of $p$ and $q$ in terms of local integrals of $p$ and $q$. 

\null

For the purpose of the proof, we  assume that the following two lemmas are true (their proof are given, respectively, in Appendices \ref{ab} and \ref{ac}).

\begin{lemma}
\label{lem1}
There exist some constants  $\lambda _0>0$,  $s_0>0$ and $C_0>0$ such that for all $\lambda \ge \lambda _0$, all $s\ge s_0$ and all 
$p\in C([0,T];L^2(\Omega ))$ with $p_t + \Delta p\in L^2 (0,T;L^2(\Omega ))$, the following  holds
\begin{multline}
\int_0^T\!\!\!\int_{\Omega } [  (s\theta )^{-1}  ( |\Delta p |^2  +|p_t|^2 ) 
+ \lambda ^2 (s\theta ) |\nabla p|^2    +  \lambda ^4 (s\theta ) ^3 |p|^2  ]e^{-2s\varphi} dxdt \\
 \le C_0 \left( \int_0^T\!\!\!\int_{\Omega } |p_t + \Delta p |^2 e^{-2s\varphi} dxdt + \int_0^T \!\!\! \int_{X(\omega  _1, t,0) }  \lambda ^4  (s\theta )^3 |p|^2 e^{-2s\varphi} dxdt \right).  \label{E99} 
\end{multline}
\end{lemma}
\noindent

\begin{lemma}
\label{lem2}
There exist some numbers $\lambda _1\ge \lambda _0$,  $s_1\ge s_0$ and $C_1>0$ such that for all $\lambda \ge \lambda _1$, all $s\ge s_1$ and
all $q\in H^1(0,T;L^2(\Omega ))$, the following  holds
\be
\int_0^T\!\!\!\int_{\Omega} ( \lambda ^2  s\theta ) |  q |^2 e^{-2s\varphi }dxdt 
\le C_1\left( 
\int_0^T\!\!\!\int_{\Omega} |q_t|^2 e^{-2s\varphi} dxdt + \int_0^T\!\!\!\int_{ X(\omega ,t,0)  } \lambda ^2 ( s\theta )^2  |q|^2 
e^{-2s\varphi} dxdt  
\right). \label{E50}
\ee  
\end{lemma}
\noindent

Applying the Carleman inequality given in Lemma  
\ref{lem1} to the heat-like equation (\ref{B1000-1}), we obtain

\begin{multline}
\int_0^T\!\!\!\int_{\Omega } [  (s\theta )^{-1}  ( |\Delta p |^2  +|p_t|^2 ) 
+ \lambda ^2 (s\theta ) |\nabla p|^2    +  \lambda ^4 (s\theta ) ^3 |p|^2  ]e^{-2s\varphi} dxdt \\
 \le C_0 \left( \int_0^T\!\!\!\int_{\Omega } |(b(x - 1 )(p - q)  |^2 e^{-2s\varphi} dxdt + \int_0^T \!\!\! \int_{X(\omega _1, t,0) }  \lambda ^4  (s\theta )^3 |p|^2 e^{-2s\varphi} dxdt \right).  \label{G1-1} 
\end{multline}
Next, we apply the Carleman inequality given by Lemma \ref{lem2} to the  ODE (\ref{B2000-1}) and obtain
\begin{multline}
\int_0^T\!\!\!\int_{\Omega} ( \lambda ^2  s\theta ) |  q |^2 e^{-2s\varphi }dxdt \\
\le C_1\left( 
\int_0^T\!\!\!\int_{\Omega} | q -p  |^2 e^{-2s\varphi} dxdt + \int_0^T\!\!\!\int_{ X(\omega ,t,0)  } \lambda ^2 ( s\theta )^2  |q|^2
e^{-2s\varphi} dxdt  
\right). \label{G2-1}
\end{multline}

Adding (\ref{G1-1}) and (\ref{G2-1}), it is not difficult to see that 

\begin{multline}
\int_0^T\!\!\!\int_{\Omega } [  (s\theta )^{-1}  ( |\Delta p |^2  +|p_t|^2 ) 
+ \lambda ^2 (s\theta ) |\nabla p|^2    +  \lambda ^4 (s\theta ) ^3 |p|^2  ]e^{-2s\varphi} dxdt + \int_0^T\!\!\!\int_{\Omega} ( \lambda ^2  s\theta ) |  q |^2 e^{-2s\varphi }dxdt  \\
 \le C \left(\int_0^T\!\!\!\int_{ X(\omega ,t,0)  } \lambda ^2 ( s\theta )^2  |q|^2 e^{-2s\varphi} dxdt  
 + \int_0^T \!\!\! \int_{X(\omega _1, t,0) }  \lambda ^4  (s\theta )^3 |p|^2 e^{-2s\varphi} dxdt \right) \label{G3-1} 
\end{multline}
for  appropriate $s \geq s_2 \ge s_1$ and $ \lambda \geq \lambda_2\ge \lambda _1$.

\null

\textit{Step 2. } {\sc Arrangements and conclusion.}\\

In this step we estimate the local integral of $p$ appearing in \eqref{G3-1} by a local integral of $q$ and some small order terms.  Finally, using semigroup theory, we finish the proof of Proposition \ref{prop1}.


The main result of this step is the following.
\begin{lemma}
\label{lem3}
There exist some numbers $\lambda _2\ge \lambda _1$, $s_2\ge s_1$ and $C_2>0$ such that for all $\lambda \ge \lambda _2$,
all  $s\ge s_2$ and all $(p_0,q_0)\in L^2(\Omega ) ^2$, the corresponding solution $(p,q)$ of system \eqref{B1000-1}-\eqref{B5000-1} fulfills
\begin{multline}
\int_0^T\!\!\!\int_{\Omega } [  (s\theta )^{-1}  ( |\Delta p |^2  +|p_t|^2 ) 
+ \lambda ^2 (s\theta ) |\nabla p|^2    +  \lambda ^4 (s\theta ) ^3 |p|^2  ]e^{-2s\varphi} dxdt + \int_0^T\!\!\!\int_{\Omega} ( \lambda ^2  s\theta ) |  q |^2 e^{-2s\varphi }dxdt  \\
 \le C_2\int_0^T\!\!\!\int_{ X(\omega ,t,0)  } \lambda ^8 ( s\theta )^7 e^{-2s\varphi} |q|^2 dxdt  .\label{G3-2} 
\end{multline}
\end{lemma}
\noindent
{\em Proof of Lemma \ref{lem3}.}  In order to prove Lemma \ref{lem3}, we just need to estimate the  $p$ appearing in the right-hand side of (\ref{G3-1} ).  For that, we 
introduce the function 
\begin{equation}\label{zeta}
\zeta (x,t) := \xi (X(x,0,t) ),
\end{equation}
where $\xi$ is a cut-off function satisfying 
\ba
&&  \xi\in C_0^\infty (\omega  ),\label{AW1}\\
&&0\le \xi (x)\le 1, \quad x\in \R ^N, \label{AW2}\\
&&\xi (x)=1,\quad  x\in \omega  _1 . \label{AW3}
\ea

We have that

\begin{equation}
\int_0^T \!\!\! \int_{X(\omega_1, t,0) }  \lambda ^4  (s\theta )^3 |p|^2 e^{-2s\varphi} dxdt  
\leq \int_0^T \!\!\! \int_{\Omega } \zeta  \lambda ^4  (s\theta )^3 |p|^2 e^{-2s\varphi}  dxdt
\label{WXYZ0}
\end{equation}

and we use (\ref{B2000-1}) to write

\begin{eqnarray}
\int_0^T \!\!\! \int_{\Omega } \zeta \lambda ^4  (s\theta )^3 |p|^2 e^{-2s\varphi} dxdt 
&=&  \int_0^T \!\!\! \int_{\Omega } \zeta  \lambda ^4  (s\theta )^3  pq e^{-2s\varphi}  dxdt  \nonumber \\
 &&\ \  +\int_0^T \!\!\! \int_{\Omega } \zeta  \lambda ^4  (s\theta )^3  p(-q_t) e^{-2s\varphi} dxdt   =: M_1 + M_2.
\end{eqnarray}

It remains to  estimate $M_1$ and $M_2$. Using Cauchy-Schwarz inequality and \eqref{AW1}-\eqref{AW2}, we have, for every $\varepsilon >0$,
\begin{equation}
|M_1| \le \varepsilon \int_0^T \!\!\! \int_{\Omega } \lambda ^4  (s\theta )^3  |p|^2 e^{-2s\varphi} dxdt  
+\frac{1}{4\varepsilon} \int_0^T \!\!\! \int_{X(\omega , t,0)} \lambda ^4  (s\theta )^3  |q|^2e^{-2s\varphi}  dxdt .\label{WXYZ1}
\end{equation}
On the other hand, integrating by parts with respect to  $t$ in $M_2$ yields

\begin{eqnarray*}
M_2 &=&\int_0^T \!\!\! \int_{\Omega} \zeta \lambda ^4  (s\theta)^3  p_tq e^{-2s\varphi}  dxdt   + 
 \int_0^T \!\!\! \int_{\Omega} \zeta \lambda ^4  (3 s^3\theta^2 \theta_t -2 s^4 \varphi _t \theta ^3) pq e^{-2s\varphi}  dxdt  \\
&& \qquad - \int_0^T \!\!\! \int_{\Omega}   \nabla\xi  (X(x,0,t) ) \cdot 
  \big( \frac{\partial X}{\partial x} \big) ^{-1} (X(x,0,t),t,0) f(x,t)  \lambda ^4  (s\theta ) ^3  pq e^{-2s\varphi}  dxdt\\
  &=:&M_2^1+M_2^2-M_2^3. 
\end{eqnarray*}
For $M_2^1$, we notice that for every $\varepsilon >0$,
\begin{equation}
|M_2^1| \le  \varepsilon \int_0^T \!\!\! \int_{\Omega } (s\theta )^{-1}  |p_t|^2 e^{-2s\varphi} dxdt  
+\frac{1}{4\varepsilon} \int_0^T \!\!\! \int_{X(\omega , t , 0)} \lambda ^8  (s\theta )^7  |q|^2e^{-2s\varphi}  dxdt .\label{WXYZ2}
\end{equation}
Since $|\theta_t| +|\varphi_t| \leq C\lambda \theta  ^2$, we infer that
 \begin{eqnarray}
|M_2^2| &\le&  C \int_0^T \!\!\! \int_{\Omega } \zeta   s^4(\lambda \theta )^5  |pq| e^{-2s\varphi} dxdt \nonumber  \\
&\le& 
\varepsilon \int_0^T \!\!\! \int_{\Omega }  \lambda ^4  (s\theta )^{3}  |p|^2 e^{-2s\varphi} dxdt  
+\frac{C}{\varepsilon s^2} \int_0^T \!\!\! \int_{X(\omega , t , 0)} \lambda ^6  (s\theta )^7  |q|^2e^{-2s\varphi}  dxdt .\label{WXYZ3}
\end{eqnarray}
Finally, $M_2^3$ is estimated as $M_1$:
\begin{equation}
|M_2^3|  
\leq  \epsilon \int_0^T\!\!\!\int_{\Omega } \lambda ^4 (s\theta ) ^3 |p|^2  e^{-2s\varphi} dxdt 
+ \frac{C}{\varepsilon} \int_0^T \!\!\! \int_{X(\omega , t,0)}  \lambda ^4(s\theta )^3 |q|^2 e^{-2s\varphi} dxdt.
\label{WXYZ4}
\end{equation}

Gathering together \eqref{G3-1} and  \eqref{WXYZ0}-\eqref{WXYZ4}  and taking $\epsilon$ small enough, we obtain \eqref{G3-2}.

Now we finish the proof of the observability inequality  \eqref{observability}.

Pick any $(p_0,q_0)\in L^2(\Omega )^2$, and denote by $(p,q)$ the solution of 
\eqref{B1000-1}-\eqref{B5000-1}. Note that $p\in C([0,T];L^2(\Omega )) \cap L^2(0,T;H^1_0(\Omega ))$ and that $q\in H^1(0,T;L^2(\Omega ))$. 
Using classical semigroup estimates, one derives at once \eqref{observability} from \eqref{G3-2}.

\qed

\section{Final comments}\label{finalcommets}

\begin{itemize}
\item {\em Another  decomposition} 

As commented in the introduction, there is another splitting of the operator ${\mathcal L} = \partial _t ^2 -\Delta -\Delta \partial _t +\partial _t$, given by 
\[
{\mathcal L} = (\partial _t -\Delta ) (\partial _ t +Id).
\]
Thus, letting 
\[
v(x,t) = y(x,t) + y_t(x,t), 
\]
we see that \eqref{A1}  may be written as
\ba
v_t -\Delta v &=& 1_{\omega(t)} h + (1-b(x))(v-y) , \label{A11}\\
y_t+y &=& v,\label{A12}
\ea
which is a coupled system of a parabolic equation \eqref{A11} and an ODE \eqref{A12}.  

This splitting was used to prove Theorem \ref{thm1} with less assumptions on the trajectories (see Section\ref{uncoupled}).

The control term $h$ acts directly in the heat equation and indirectly in the ODE through the coupling term $v$. The problem can be treated direct�y as such, with requires further work at the level of the dual observability problem since both components of the adjoint system will be needed to be observed by partial measurements only on one of its components. The problem can also be addressed 
 incorporating in \eqref{A12} an additional auxiliary control  acting directly in the ODE. This leads to the system
\ba
v_t -\Delta v &=& 1_{\omega(t)}  h + (1-b(x))(v-y) \label{A13-1}\\
y_t+y &=& 1_{\omega(t)} k +  v,\label{A14}
\ea
where $(v,y)\in L^2(\Omega )^2$ is the state function to be controlled, and $(h,k)\in L^2(0,T;L^2(\Omega )^2)$ is the control input. 

Once the controllability of this system is proved, when going back to the original viscoelasticity equation, one gets
\be
\label{A40}
y_{tt}-\Delta y - \Delta y_t +b(x)y_t= 1_{\omega(t)} \, [h - (1-b(x)) k]  +(\partial _t-\Delta)[1_{\omega(t)} k]. 
\ee
But, then, the second control $1_{\omega(t)} k$ enters under the action of the heat operator. It is then necessary to ensure that the control $k$ is smooth enough and, furthermore, to replace in \eqref{A14}  the cut-off function $1_{\omega(t)}$ by a regularized version. 
These are technicalities that can be overcame with further work. To be more precise, 
the control in   \eqref{A40} takes the form ${\bf 1}_{X(\omega ,t,0)} (x)  \tilde h$ with $\tilde h\in L^2(0,T;L^2(\Omega ))$, provided that both $ h, \, k \in L^2(0,T;L^2(\Omega ))$ and
\[
k\in H^1(0,T; L^2(\Omega ))\cap L^2(0,T;H^2(\Omega)).
\] 
Therefore special attention has to be paid to obtain smooth controls for the transport equation (see Section \ref{uncoupled}).

\item {\em Manifolds without boundary} 

The lack of propagation properties of the ODE \eqref{A14bis}  in the space variable requires the control  to move in time. As we mentioned in the introduction, 
through a suitable change of variables, this is equivalent to keeping the support of the control fixed but replacing the ODE by a transport equation. Obviously, attention 
has to be paid to the Dirichlet boundary conditions when performing this change of variables. Of course, this is no longer an issue when the model is considered 
in a smooth manifold without boundary. As an example of such a situation we consider  the periodic case in the torus
\be
\label{TN}
x\in \T ^N :=\R ^N / \Z ^N.
\ee
For a moving control with a constant velocity 
$\omega(t)= \{ x-ct; \  x\in \omega \}$, $c\in \R ^N\setminus\{0\}$, system \eqref{A13bis}-\eqref{A14bis}  can be put in the form of a coupled system of parabolic-hyperbolic equations
\begin{eqnarray}
v_t -\Delta v - c\cdot \nabla v + (b(x+ct)-1)v&=& w  \label{A15}\\
w_t- c\cdot \nabla w  + w &=& 1_\omega(x)\tilde h +  (b(x+ct)-1) v   \label{A16}
\end{eqnarray}
by letting 
\ba
v(x,t) &=& y(x+ct, t), \label{A17a} \\
w(x,t) &=& z(x+ct,t), \label{A17b} \\
\tilde h(x,t) &=& h(x+ct,t). \label{A17c} \\
\ea
The system is now constituted by the coupling between a heat and a transport equation with control $\tilde h $  
with fixed support. Once more, the problem now can be treated by means of the classical duality principle between the controllability problem and the observability property of the adjoint system.  The later was solved  in \cite{MRR} in $1-d$ using Fourier analysis techniques and in this paper we do it using Carleman inequalities.

Note that the Carleman approach developed in this paper cannot be applied as it is to the periodic case. 
Consider for instance the case of the torus $\T$. A weight $\psi \in C^\infty (\T \times (0,T))$ as in Lemma \ref{weight} 
does not exist, because of the periodicity in $x$ (see Figure \ref{fig4}.) However, it is well known that the periodic case can be deduced from both the Dirichlet case and the 
Neumann case (using classical extensions by reflection, see e.g. \cite{RZ2010}). Even if the Neumann case was not considered in this paper, it is likely that it could be treated
 in much the same way as we did for the Dirichlet case. 


\end{itemize}

\section{Open problems and further questions}\label{openproblems}

The main result of this paper concerns the controllability of a coupled system consisting on a heat equation and an ODE. By addressing the dual problem of observability and making the controller/observer move in time, this ends being very close to the problem of observability of a coupled system of a heat equation and a first order transport equation. The techniques we have developed here are inspired in the work \cite{AT} where the key point was to use the same weight function for the Carleman inequality in both the heat and the transport equation. 

The system under consideration, coupling a heat and a hyperbolic equation, is close to that of thermoelasticity that was considered in \cite{LZ}. But, there, the problem was only dealt with in the case of manifolds without boundary, by means of spectral decomposition techniques allowing to decouple the system into the parabolic and the hyperbolic components. As far as we know, a complete analysis of the system of thermoelasticity using Carleman type inequalities seems to be not developed so far.

The structure of the parabolic-transport system we consider is also, in some sense, similar to the one considered in \cite{EGGP} for the $1-d$ compressible Navier-Stokes equation although, in the latter, the system is of nonlinear nature requiring significant extra analysis beyond the linearized model.

Our analysis is also related to recent works on the control of parabolic equations with memory terms as for instance in \cite{GI}. Note that the system \eqref{A13bis}-\eqref{A14bis} in the particular case $b \equiv 1$ and $z(0)\equiv 0$, in the absence of the control $h$ and the addition of a  control of the form $1_{\omega(t)} k$ in the first equation,  can be written as an integro-differential equation
\ba
y_t -\Delta y  + (b-1) [y - \int_0^t e^{s-t} y(x, s) ds] = 1_{\omega(t)} k \label{A13x}
\ea
This system is closely related to the one considered in \cite{GI}. There it is shown that the system lacks to be null controllable. This is in agreement with our results that, in the particular case under consideration, show also that a moving control could bypass this limitation. It would be interesting to analyze to which extent this idea of controlling by moving the support of the control can be of use for more general parabolic equations with memory terms.

In this paper we have shown the null controllability of a linear system which consists of a  parabolic equation and an ordinary differential equation that 
arise from the identification of the parabolic and hyperbolic parts of system (\ref{A1})-(\ref{A1bis}). Besides, coupled systems consisting of parabolic equations and ode's are important since they appear in biological models of chemotaxis or interactions between cellular process and diffusing growth factors (see \cite{Hors}, \cite{MC-K-S},  \cite{Rascle-Ziti} and references therein). Systems governing these phenomena are, in general, non linear and have the form 

\begin{eqnarray}
u_t  &=& f(u,v), \label{x1}\\
v_t  &=& D\Delta v + g(u,v),  \label{x2}
\end{eqnarray}
where $v$ and $u$ are vectors, $D$ is a diagonal matrix with positive coefficients and $f$ and $g$ are real functions.

Other area where coupled  parabolic-ode systems play a major role is electrocardiology (see \cite{B-K}, \cite{CF-P}, \cite{Z-J} and references therein). Here the cardiac activity is described by the bidomain model, which consists  of a system of two degenerate parabolic reaction-diffusion equations, representing the intra and extracellular potential in the cardiac muscle, coupled with a system of ordinary differential equations representing the ionic currents flowing through the cellular membrane. The bidomain model is given by 
\begin{eqnarray}
\chi C_m v_t -\textrm{div}\, (D_i \nabla u_i) + \chi I_{ion}(v,w) &=& I^i_{app},  \label{x3}\\
-\chi C_m v_t -\textrm{div}\, (D_e \nabla u_e) - \chi I_{ion}(v,w) &=& -I^e_{app},  \label{x4}\\
w_t -R(v,w) = 0,  \label{x5}
\end{eqnarray}
where $u_i$ and $u_e$ are the intra and extracellular potentials, $v$ is the transmembrane potential, $\chi$ is the ratio  of membrane area per tissue volume, $C_m$ is the surface capacitance of the membrane, $I_{ion}$ is the ionic current,   $I^{i,e}_{app}$ is an applied current and $D_{i,e}$ are conductivity tensors.

Concerning controllability of coupled parabolic-ode systems, just a few results for some particular systems are known (see \cite{D-FC} and \cite{Liu-T-T}  for the controllability of a simplified one-dimensional model for the motion of a rigid body in a viscous fluid). We believe that ideas presented in this paper can be used for the study of the controllability for other systems of parabolic-ode equations, such as (\ref{x1})-(\ref{x2}) and (\ref{x3})-(\ref{x5}).

\appendix 

\section{Proof of Lemma \ref{weight}}\label{aa}
\begin{proof}
Pick any $\delta <\min(t_1,T -t_2,T/2)$. 
We search $\psi$ (see Figure \ref{fig4}) in the form 
\be
\label{E100}
\psi (x,t) =\psi_1(x,t) + C_2\psi _2(x,t) + C_3 
\ee
where, roughly, $\psi_1$ fulfills \eqref{P1}, $\psi_2$ fulfills \eqref{P2}-\eqref{P4} together with $\nabla \psi_2\equiv 0$ outside $X(\omega_1,t,0)$, and $C_2,C_3$ are 
(large enough) positive constants.\\
\textit {Step 1.} {\sc  Construction of $\psi _1$.}\\
Let $\Gamma \in C^\infty ([0,T];\R ^N)$ be as in \eqref{A3a}, and let $\varepsilon >0$ be such that 
\[
B(\Gamma (t),3\varepsilon ) \subset X(\omega _0, t,0)\cap \Omega, \qquad t\in [0,T]. 
\]
We choose a vector field $\tilde f\in C^\infty (\R ^N\times [0,T];\R ^N)$ such that 
\[
\tilde f(x,t) =
\left\{ 
\begin{array}{ll}
\dot\Gamma (t) \quad & \text{ if } t\in [0,T], \ x\in B(\Gamma (t), \varepsilon ), \\ 
0                                     & \text{ if } t\in [0,T], \ x\in \R ^N \setminus B(\Gamma (t), 2\varepsilon ).  
\end{array}
\right. 
\]
Let $\tilde X$ denote the flow associated with $\tilde f$; that is, $\tilde X$ solves
\begin{eqnarray*}
\frac{\partial \tilde X}{\partial t} (x,t,t_0) = \tilde f (\tilde X (x,t,t_0),t),&& \quad (x,t,t_0)\in \R ^N\times [0,T]^2,\\
{\tilde X}(x,t_0,t_0) = x,&& \quad  (x,t_0) \in \R ^N \times [0,T].
\end{eqnarray*}
Note that 
\begin{eqnarray*}
\tilde X (y+\Gamma (0),t,0) = y+\Gamma (t)  &\text{if}& \ (y,t)\in B(0,\varepsilon ) \times [0,T],\\
\tilde X(x,t,t_0) =x &\text{if}& \ \text{dist } (x,\partial \Omega ) <\varepsilon, \ (t,t_0) \in [0,T]^2.  
\end{eqnarray*}
By a well-known result (see \cite[Lemma 1.2]{Imanuvilov}), there exists a function $\tilde \psi \in C^\infty(\overline{\Omega })$ such that 
\begin{eqnarray*}
\tilde\psi (x ) >0 &\text{if}& \ x\in \Omega;\\
\tilde\psi (x) =0 &\text{if}&\ x\in\partial \Omega ;\\
\nabla \tilde \psi (x)\ne 0 &\text{if}& \ x\in \overline{\Omega} \setminus B(\Gamma (0),\varepsilon ).    
\end{eqnarray*}
Actually, the function $\tilde\psi$ given in \cite{Imanuvilov} is only of class $C^2$, but the regularity $C^\infty$ can be obtained by mollification
with a partition of unity
(see e.g. \cite[Lemma 4.2]{RZ2009}). 
Let us set 
\[
\psi _1(x,t) = \tilde \psi (\tilde X (x,0,t)).
\]
Then $\psi _1\in C^\infty(\overline{\Omega} \times [0,T])$ and it fulfills
\ba
\psi _1(x,t) >0 &\text{if}&\ (x,t)\in \Omega \times [0,T], \label{FF1}\\
\psi _1 (x,t)=0 &\text{if}& \ (x,t)\in \partial\Omega \times [0,T],\label{FF2}\\
\nabla \psi _1 (x,t) = \nabla \tilde \psi (\tilde X(x,0,t)) \frac{\partial \tilde X}{\partial x} (x,0,t) \ne 0 &\text{if}& \ x\in 
\overline{\Omega } \setminus X(\omega _0,t,0).\label{F3} 
\ea
For \eqref{F3}, we notice that if we write $x=\tilde X (\tilde x,t,0)$, then $\tilde x=\tilde X(x,0,t)$ hence
\[
\nabla \tilde \psi (\tilde X (x,0,t)) =\nabla \tilde \psi (\tilde x )\ne 0
\]
if $\tilde x \not\in B(\Gamma (0),\varepsilon )$, which is equivalent to $x\not\in B(\Gamma (t),\varepsilon )$. The last condition is satisfied when $x\not\in X(\omega _0,t,0)$.\\
\textit {Step 2.} {\sc Construction of $\psi _2$.}\\
From \eqref{A3c}, \eqref{A3d}, and \eqref{A3e}, we can pick two curves $\gamma _1\in C^0( [0,t_2) ;\Omega )$ and $\gamma _2 \in C^0( (t_1,T] ;\Omega )$ such that 
\begin{eqnarray*}
\gamma _1 (t) \not \in \overline{X(\omega _0,t,0)},&& \qquad 0\le t < t_2,\\
\gamma _2 (t) \not\in \overline{X(\omega _0,t,0)},&& \qquad t_1 < t\le T. 
\end{eqnarray*}  
We infer from \eqref{A3e} that for any $t\in (t_1,t_2)$, $\gamma _1 (t)$ and $\gamma _2(t)$ do not belong to the same connected component of
$\Omega \setminus \overline{X(\omega _0,t,0)}$. Let $\Omega _1(t)$ (resp. $\Omega _2(t)$) denote the connected component of $\gamma _1(t)$ 
(resp. $\gamma _2(t)$) for $0\le t < t_2$ (resp. for $t_1 < t\le T$). Clearly 
\[
\Omega \setminus \overline{X(\omega _0, t,0)} 
=\left\{ 
\begin{array}{ll}
\Omega _1 (t),\quad &\text{if} \ 0\le t\le t_1,\\
\Omega _1(t)\cup \Omega _2(t), \quad &\text{if} \ t_1<t <t_2,\\
\Omega _2 (t), \quad &\text{if}\  t_2\le t\le T.
\end{array}
\right.
\]
Set $\Omega _1(t)=\emptyset$ for $t\in [t_2,T]$, and $\Omega _2(t)=\emptyset$ for $t\in [0,t_1]$. 
Let  $\psi _2\in C^\infty (\overline{\Omega}\times [0,T])$ with 
\begin{eqnarray*}
\psi _2(x,t) &=& t \left( {\bf 1}_{\Omega _1 (t)} (x) - {\bf 1}_{\Omega _2 (t) } (x)\right) \quad \text{ for } 
t\in [0,T], \ x\in \Omega \setminus X(\omega _1,t,0),\\
\frac{\partial \psi_2}{\partial n} &=&0 \quad \text{ for } (x,t)\in \partial \Omega \times [0,T].
\end{eqnarray*}
Such a function $\psi _2$ exists, since by \eqref{A3-x}
\[
\inf_{t_1<t<t_2} \textrm{dist}\ \big(\Omega _1(t) \setminus X(\omega _1, t, 0 ) , \Omega _2(t) \setminus X(\omega _1, t, 0 )\big) >0.
\]
Then
\[
\frac{\partial \psi _2}{\partial t}
=\left\{ 
\begin{array}{ll}
1\quad &\text{if}\  \ 0\le t < t_2 \ \text{ and }\    x\in \Omega _1 (t) \setminus X(\omega _ 1, t, 0) ,\\
-1\quad &\text{if} \ \ t_1 <  t\le T \ \text{ and }  \ x\in \Omega _2 (t) \setminus X(\omega _1 , t, 0) .
\end{array}
\right.
\]
and 
\[
\nabla \psi _2(x,t) =0 \quad \text{ if }\  x\in \Omega \setminus X (\omega _1,t,0) .
\]
Note that \eqref{P2}-\eqref{P4} are satisfied for $\psi _2$. Note also that for any pair $(\tau _1, \tau _2)$ with $0\le \tau _1 < \tau _2 \le T$, the set
\[ K_{\tau _1,\tau _2} := \{ (x,t)\in\R ^{N+1}; \ \ \tau _1\le t\le \tau _2, \  x\in \overline{\Omega} \setminus X(\omega _1 , t, 0) \} \]    
is {\em compact}.  
 
\noindent
\textit {Step 3.} {\sc Construction of $\psi$.} \\
Let $\psi$ be defined as in \eqref{E100}, with $C_2>0$ and $C_3$ to be determined. Then \eqref{P1} and \eqref{P5} are satisfied. We pick $C_2$
large enough for \eqref{P2}-\eqref{P4} to be satisfied. Finally, \eqref{P6} is satisfied for $C_3$ large enough. 
\end{proof}

\section{Proof of Lemma \ref{lem1} } \label{ab}

{\em Proof of Lemma \ref{lem1}.}  The method of the proof is widely inspired from \cite{FI}, and the computations are presented as in 
\cite[Proof of Proposition 4.3]{RZ2009}. 

Let $v=e^{-s\varphi}p$ and $P=\partial _t + \Delta $. Then 
\[
e^{-s\varphi} Pp  = e^{-s\varphi } P(e^{s\varphi } v ) = P_s v + P_av
\]
where 
\begin{eqnarray}
P_s v &=&   \Delta v + (s\varphi _t + s^2 |\nabla \varphi |^2 )v , \label{P50} \\
P_a v &=&  v_t + 2s \nabla \varphi  \cdot  \nabla v  + s (\Delta  \varphi ) v \label{P60}
\end{eqnarray}
denote the (formal) self-adjoint and skew-adjoint parts of $e^{-s\varphi } P(e^{s\varphi} \cdot )$, respectively.  
It follows that 
\be
\label{B1500}
||e^{-s\varphi} Pp||^2 = ||P_s v||^2 + ||P_a v||^2 +  2(P_s v,P_a v)
\ee
where $(f,g)=\int_0^T\!\!\!\int_{\Omega } fg \,dxdt$, $|| f ||^2 = (f,f)$. In the sequel, $\noindent \int_0^T\!\!\!\int_{\Omega }f(x,t) dxdt $ is denoted $ \noindent \ii f$ for the sake of shortness. We have
\begin{multline}
(P_s v,P_a v) = \big(  \Delta v , v_t \big) + \big(  \Delta v, 2s \nabla \varphi \cdot \nabla v \big) + 
(\Delta v , s (\Delta \varphi  ) v) 
+\big( s\varphi _t v +  s^2|\nabla \varphi |^2 v, v_t) \\
 +  (s\varphi _t  v + s^2 |\nabla \varphi |^2 v, 2s \nabla \varphi \cdot \nabla v)   +  (s\varphi _t v +  s^2 |\nabla \varphi |^2 v, s (\Delta \varphi  ) v)
= : I_1 + I_2 + I_3 + I_4 + I_5  + I_6. \label{B0}
\end{multline}
First, observe that 
\be
\label{B1}
I_1 = -\ii \nabla v\cdot \nabla v_t =0. 
\ee
Using the convention of repeated indices and denoting $\partial _i=\partial /\partial x_i$, we obtain that 
\begin{eqnarray*}
I_2&= &2s \ii \partial _j^2 v\, \partial _i \varphi\,  \partial _i v \\
&=& -2s \ii \partial _j v (\partial _j\partial _i \varphi \partial _i v    +    \partial _i \varphi \partial _j \partial _i v  ) 
+2s \int_0^T\!\!\!\int_{\partial \Omega} (\partial _j v) n_j \partial _i \varphi \partial _i v d\sigma .
\end{eqnarray*}
Since $v=0$ for $(x,t)\in \partial \Omega \times (0,T)$, $\nabla v=(\partial v/\partial n) n$, so that $\nabla \varphi \cdot  \nabla v = (\partial \varphi  / \partial n) (\partial v / \partial n) $ and 
\[
\int_0^T\!\!\!\int_{\partial \Omega} (\partial _j v) n_j \partial _i \varphi \partial _i v \, d\sigma 
= \int_0^T\!\!\! \int_{\partial \Omega } (\partial \varphi/\partial n) | \partial v/ \partial n|^2 d\sigma . 
\]
It follows that 
\ba
I_2&=&-2s \ii \partial _j \partial _i \varphi \partial _j v\partial _i v  - s\ii \partial _i \varphi \partial _i (|\partial _j v|^2)    
+2s  \int_0^T\!\!\! \int_{\partial \Omega } (\partial \varphi/\partial n) | \partial v/ \partial n|^2 d\sigma 
\nonumber \\
 &=&-2s \ii \partial _j \partial _i \varphi \partial _j v\partial _i v  + s\ii \Delta \varphi   |\nabla v|^2 
  +s  \int_0^T\!\!\! \int_{\partial \Omega } (\partial \varphi/\partial n) | \partial v/ \partial n|^2  d\sigma   \label{B2}
\ea
On the other hand, integrations by parts in $x$ yields 
\be
I_3 = -s\ii \nabla v \cdot \big(v \nabla (\Delta \varphi)  + (\Delta \varphi )\nabla v \big) = s \ii \Delta ^2 \varphi \frac{|v|^2}{2}  -s\ii \Delta \varphi |\nabla  v|^2
\label{B3}
\ee
and integration by parts with respect to $t$ gives
\[
I_4 = -\ii (s\varphi _{tt} +s^2 \partial _t |\nabla \varphi | ^2) \frac{|v|^2}{2}\cdot \label{B4}
\]
Integrating by parts with respect to $x$ in $I_5$ yields 
\be
I_5= -\ii s^2 \nabla \cdot (\varphi _t \nabla \varphi) |v|^2 -\ii s^3 \nabla \cdot (|\nabla \varphi |^2\nabla \varphi  ) |v|^2. 
\label{B5}
\ee
Gathering \eqref{B0}-\eqref{B5}, we infer that 
\begin{multline*}
2(P_sv,P_av) =-4s \ii \partial_j\partial_i \varphi \partial _j v\partial _i v  
+2s  \int_0^T\!\!\! \int_{\partial \Omega } (\partial \varphi/\partial n) | \partial v/ \partial n|^2  d\sigma \\
+\ii |v|^2 [ s(\Delta ^2 \varphi - \varphi _{tt}) -2s^2 \partial _t |\nabla \varphi |^2 -2s^3 \nabla \varphi \cdot \nabla 
|\nabla \varphi |^2 ].
\end{multline*}

Consequently, \eqref{B1500} may be rewritten 
\begin{multline*}
||e^{-s\varphi } Pp||^2 = ||P_s v||^2 + ||P_a v|| ^2 -4s \ii \partial_j\partial_i \varphi \partial _j v\partial _i v
 + 2s  \int_0^T\!\!\! \int_{\partial \Omega } (\partial \varphi/\partial n) | \partial v/ \partial n|^2  d\sigma  \\
+\ii |v|^2 [ s(\Delta ^2 \varphi - \varphi _{tt}) -2s^2 \partial _t |\nabla \varphi |^2 -2s^3 \nabla \varphi \cdot \nabla 
|\nabla \varphi |^2 ].
\end{multline*}

{\sc Claim 1.} There exist some numbers $\lambda _1 >0$, $s_1>0$ and  $A\in (0,1)$ such that for all $\lambda \ge \lambda _1$ and all $s\ge s_1$, 
 \begin{multline}
\ii |v|^2 [ s(\Delta ^2 \varphi - \varphi _{tt}) -2s^2 \partial _t |\nabla \varphi |^2  -2s^3 \nabla \varphi \cdot \nabla 
|\nabla \varphi |^2 ] \\
+A^{-1}\lambda s^3 \int_0^T\!\!\! \int_{X( \omega  _1,t,0) } (\lambda \theta ) ^3 |v|^2 \ge A\lambda s^3 \ii (\lambda \theta ) ^3 |v|^2. 
\label{B10} 
\end{multline}

\noindent
{\em Proof of Claim 1.} 
Easy computations show that 
\be
\label{B400}
\partial _i \varphi =-\lambda g(t) e^{\lambda \psi} \partial _i \psi, \qquad  \partial_j\partial _i\varphi 
=-g(t) e^{\lambda \psi} (\lambda ^2 \partial _i \psi \partial _j \psi + \lambda \partial _j\partial _i \psi)  
\ee
and 
\[
-\nabla |\nabla \varphi |^2 \cdot \nabla \varphi = -2 (\partial_j\partial _i \varphi) \partial _i\varphi\partial_j \varphi
=2(\lambda g e^{\lambda \psi})^3 (\lambda |\nabla \psi |^4 + \partial _j\partial _i \psi \partial _i \psi  \partial _j \psi). 
\]
It follows from  \eqref{P1} that  for $\lambda $ large enough, say $\lambda \ge \lambda _1$, we have that
\ba
-\nabla |\nabla \varphi |^2 \cdot \nabla \varphi &\ge& A \lambda (\lambda \theta )^3, \qquad t\in [0,T],
\  x\in \overline{\Omega}\setminus X ( \omega _1,t,0)   \\
|\nabla |\nabla \varphi |^2 \cdot \nabla \varphi |  &\le& A^{-1} \lambda (\lambda \theta )^3, \qquad t\in [0,T],\  x\in  X ( \omega _1,t,0)
\ea
for some constant $A\in (0,1)$. 
According to \eqref{P6}, we have for some constant $C>0$ 
\[
|\Delta ^2 \varphi |  + |\varphi _{tt}| + |\partial _t |\nabla \varphi |^2| \le C \lambda (\lambda \theta ) ^3, \qquad t\in [0,T], \ x\in \overline{\Omega}.
\]
Therefore, we infer that for $s$ large enough, say $s\ge s_1$, and for all $\lambda \ge \lambda _1$
we have that 
\begin{eqnarray*}
s(\Delta ^2 \varphi - \varphi _{tt}) -2s^2 \partial _t |\nabla \varphi |^2 -2s^3 \nabla \varphi \cdot \nabla 
|\nabla \varphi |^2   &\ge & A \lambda  s^3 ( \lambda  \theta  )^3, \quad t\in [0,T], \ x\in \overline{\Omega} \setminus X(\omega _1,t,0)\\ 
|s(\Delta ^2 \varphi - \varphi _{tt}) -2s^2 \partial _t |\nabla \varphi |^2 -2s^3 \nabla \varphi \cdot \nabla 
|\nabla \varphi |^2  |
 &\le &  3 A^{-1} \lambda  s^3 ( \lambda  \theta )^3, \quad t\in [0,T], \ x\in  X(\omega _1,t,0). 
\end{eqnarray*}
This gives \eqref{B10} with a possibly decreased value of  $A$.\qed

Thus, using the fact that $\partial \varphi/\partial n \ge 0$ on $\partial \Omega$ by \eqref{P5}, we conclude that 
\begin{multline}
||P_s v||^2 + ||P_a v||^2 + A\lambda s^3 \ii (\lambda  \theta)^3 |v|^2 \\
\le ||e^{-s\varphi} Pp||^2 + 4s \ii \partial _j\partial _i \varphi \partial _j v  \partial _i v
+A^{-1} \lambda s^3 \int_0^T\!\!\! \int_{ X( \omega _1,t,0)  } (\lambda \theta ) ^3 |v|^2. 
\label{B15}
\end{multline}
{\sc Claim 2.} There exist some numbers $\lambda _2\ge \lambda _1$ and $s_2\ge s_1$ such that for all $\lambda \ge \lambda _2$ and all $s\ge s_2$, 
\be
\label{B16}
\lambda s \ii (\lambda \theta )   |\nabla v|^2 + \lambda s^{-1} \ii (\lambda \theta )^{-1} |\Delta v|^2 \le C\left( 
s^{-1} ||P_s v||^2 + \lambda s^3 \ii (\lambda \theta )^3 |v|^2 \right) .
\ee
{\em Proof of Claim 2.} By \eqref{P50}, we have
\ba
s^{-1} \ii (\lambda \theta )^{-1} |\Delta v| ^2 
&=& s^{-1} \ii (\lambda \theta )^{-1} |P_ s v - s\varphi _t v - s^2 |\nabla \varphi | ^2 v |^2 \nonumber  \\
&\le& Cs^{-1} \ii (\lambda \theta) ^{-1} \big( 
|P_s v| ^2 + s^2 |\varphi _t|^2 |v|^2    + s^4 (\lambda \theta ) ^4 |v|^2 \big)   \nonumber \\
&\le& C\left(  \frac{||P_s v||^2 }{\lambda s}   +s^3 \ii (\lambda \theta )^3 |v|^2 \right)\label{B20} 
\ea
provided that $s$ and $\lambda $ are large enough, where we used \eqref{P6} in the last line to bound $\varphi _t$.  On the other hand, 
\ba
\lambda s \ii (\lambda \theta )  |\nabla v|^2 &=& \lambda s \{
\ii (\lambda \theta )  (-\Delta v) v - \ii (\nabla (\lambda \theta )  \cdot \nabla  v) v \}  \nonumber  \\
&\le& \frac{\lambda}{2s}\ii (\lambda \theta )^{-1} |\Delta v|^2 + \frac{\lambda s^3}{2} \ii (\lambda \theta )^3 |v|^2 +\frac{\lambda s}{2} \ii \Delta (\lambda \theta )  |v|^2 \nonumber \\
&\le& C \left( s^{-1} ||P_s v||^2 +\lambda s^3 \ii (\lambda \theta )^3 |v|^2 \right) \label{B21}
\ea
by \eqref{B20}, provided that $s\ge s_2\ge s_1$ and  $\lambda \ge \lambda _2\ge \lambda _1$. 
Then \eqref{B16} follows from \eqref{B20}-\eqref{B21}.\qed 

We infer from \eqref{B15}-\eqref{B16} that  
\begin{multline}
||P_a v||^2 +\lambda s \ii (\lambda \theta )  |\nabla v|^2  + \lambda s^{-1} \ii (\lambda \theta )^{-1} |\Delta v|^2 
+\lambda s^3 \ii (\lambda \theta )^3 |v|^2 \\
\le C \left( ||e^{-s\varphi } Pp||^2 + 4s \ii \partial _j\partial _i \varphi \partial _j v\partial _i v   
 + A^{-1} \lambda s^3 \int_0^T\!\!\! \int_{ X(\omega _1,t,0) } (\lambda \theta )^3 |v|^2\right) .\label{B24}
\end{multline}
By \eqref{B400},
\[
s\ii \partial _j \partial _i \varphi \partial _j v\partial _i v \le -s\lambda \ii g(t)e^{\lambda \psi }\partial _j\partial _i \psi \partial _j v\partial _i v 
\le C s \ii (\lambda \theta )  |\nabla v |^2. 
\]
Therefore, for $\lambda$ large enough and $s\ge s_2$, 
\begin{multline}
||P_a v||^2 + \lambda s^3 \ii (\lambda \theta )^3 |v|^2  + \lambda s \ii (\lambda \theta ) |\nabla v|^2 + 
\lambda s^{-1} \ii (\lambda \theta )^{-1} |\Delta v|^2 \\
\le C \left( ||e^{-s\varphi} Pp||^2 
+\lambda s^3 \int_0^T\!\!\!\int_{ X(\omega _1,t,0) } (\lambda \theta )^3 |v|^2 
\right) . \label{B30}
\end{multline}
Using \eqref{P60} and \eqref{B30}, we see that for $\lambda$ large enough and $s\ge s_2$ 
\begin{eqnarray*}
\lambda s^{-1} \ii (\lambda \theta )^{-1} |v_t|^2 &\le & C\lambda s^{-1} \ii (\lambda \theta )^{-1} 
\big(    |P_a v | ^2     + s^2 |\nabla \varphi |^2 | \nabla v|^2 + s^2 |\Delta \varphi |^2 |v|^2   \big)  \\
&\le& C \left(  ||e^{-s\varphi } Pp||^2 + \lambda s^3 \int_0^T\!\!\! \int_{ X(\omega _1,t,0) } (\lambda \theta )^3 |v|^2\right)    .
\end{eqnarray*}
Hence, there exists some number $\lambda _3\ge \lambda _2$ such that for all $\lambda \ge \lambda _3$ and all $s\ge s_2$, we have
\begin{multline}
\lambda s^3 \ii (\lambda \theta ) ^3 |v|^2 + \lambda s \ii (\lambda \theta )  |\nabla v|^2 + \lambda s^{-1} \ii (\lambda \theta)^{-1} (|\Delta v| ^2 + |v_t|^2) \\
\le C \left( ||e^{-s\varphi} Pp||^2 + \lambda s^3 \int_0^T\!\!\! \int_{X( \omega _1,t,0)} (\lambda \theta )^3 |v|^2\right). \label{B32}
\end{multline} 
Replacing $v$ by $e^{-s\varphi}p$ in \eqref{B32} gives
at once \eqref{E99}. The proof of Lemma \ref{lem1} is complete.  \qed

\section{Proof of Lemma \ref{lem2}}\label{ac}
{\em Proof of Lemma \ref{lem2}.} The proof is divided into three parts corresponding to the estimates for $t\in [0,\delta ]$, for $t\in [\delta ,T -\delta ]$ and for $t\in [T-\delta , T]$.  
The estimates for $t\in [0,\delta]$ and for $t\in [T-\delta, T]$ being similar, we shall prove only the first ones.

Let $v=e^{-s\varphi} q$. Then 
\be
e^{-s\varphi } q_t = e^{-s\varphi} (e^{s\varphi } v)_t 
=s\varphi _t v  + v_t  
=: P_s v + P_a v.
\ee

\noindent
{\sc Claim 3.} 
\begin{multline}
\int_0^\delta \!\!\! \int_{\Omega } \lambda ( s\theta ) ^2   | v |^2  dxdt 
\le 
C\left( \int_0^\delta\!\!\! \int_{\Omega } 
\lambda ^{-1} | e^{-s\varphi} q_t |^2 dxdt    \right. \\
\left. + 
\int_{\Omega } [(1-\zeta )^2 ( s\theta ) |v|^2]_{\vert t=\delta }dx + \int_0^\delta \!\!\!\int_{ X(\omega   ,t,0) }  \lambda ( s\theta ) ^2 
 |v|^2 dxdt   \right), \label{F0} 
\end{multline}
where $\zeta $ is the function introduced in \eqref{zeta}.

To prove the claim, we compute in several ways
\[
I:=\int_0^\delta \!\!\! \int_{\Omega } (e^{-s\varphi} q_t )(1-\zeta )^2  s\theta v \, dxdt. 
\]
We split $I$ into
\[
I = \int_0^\delta \!\!\!\int_{\Omega} (P_s v) (1-\zeta )^2  s \theta v \, dxdt + \int_0^\delta \!\!\! \int_{\Omega } (P_a v) (1-\zeta )^2  s\theta v \, dxdt =: I_1 + I_2. 
\]
Then 
\begin{eqnarray*}
I_1 &=& \int_0^\delta \!\!\! \int_{\Omega}  (1-\zeta )^2  s^2 \varphi _t  \theta v^2 \, dxdt \\
&=& \int_0^\delta \!\!\! \int_{\Omega} [g' (e ^{\frac{3}{2}\lambda  ||\psi||_{L^\infty}} -e^{ \lambda \psi}) - g  \lambda \psi _t e^{\lambda  \psi} ](1-\zeta )^2  s^2 g e^{\lambda \psi } v^2\, dxdt.  
\end{eqnarray*}
On the other hand
\begin{eqnarray*}
I_2 &=& \int_0^\delta\!\!\! \int_{\Omega} (1-\zeta ) ^2 ( s g e^{\lambda \psi} vv_t) \, dxdt\\
&=&\frac{1}{2} \int_{\Omega} [(1-\zeta ) ^2  s\theta  |v|^2 ]_{\vert t=\delta} dx  -\int_0^\delta \!\!\! \int_{\Omega}  s [g'e^{ \lambda \psi} + g  \lambda \psi _t e^{\lambda  \psi}] (1-\zeta  )^2 \frac{v^2}{2} dxdt \\
&&\qquad 
- \int_0^\delta \!\!\! \int_{\Omega}   (1-\zeta  ) \nabla \xi (X(x,0,t))   \cdot  \big( \frac{\partial X}{\partial x} \big) ^{-1} (X(x,0,t),t,0) f(x,t)   s\theta v^2 dxdt,      
\end{eqnarray*}
where we used the fact that $[\theta |v|^2]_{| t= 0 }=0$. 
Clearly, since $\theta \ge 1$, for $s\ge 1$
\begin{multline*}
| \int_0^\delta \!\!\! \int_{\Omega}   (1-\zeta  )\nabla \xi ( X(x,0,t))   \cdot 
 \big( \frac{\partial X}{\partial x} \big) ^{-1} (X(x,0,t),t,0) f(x,t) s\theta v^2 dxdt  | \\
\le C \int_0^\delta\!\!\!\int_{X(\omega ,t,0)} (s\theta )^2 |v|^2 dxdt.
\end{multline*}
On the other hand, using \eqref{P3}, we see that  there exist some constants  $C>0$  and $s_1\ge s_0$ such that  for all $s\ge s_1$ and all $\lambda \ge \lambda _0>0$, it holds
\begin{eqnarray*}
g \lambda \psi _t  e^{ \lambda \psi} ( s^2 g e^{\lambda \psi }+ \frac{ s }{2})  \ge C \lambda
 ( s\theta)^2  ,\quad && t\in (0,\delta ),\ x\in \overline{\Omega} \setminus  X(\omega _1 ,t,0) \\
-g'(t)\left( (e^{\frac{3}{2}\lambda  \Vert \psi\Vert _{L^\infty}}  - e^{\lambda \psi } )  s^2 g e^{ \lambda \psi} -\frac{ s}{2} e^{ \lambda \psi}  \right) >0\quad &&  t\in (0,\delta ), \ 
x\in \overline{\Omega} \setminus X (\omega _1 ,t,0) . 
\end{eqnarray*} 
It follows that for some positive constant $C'>C$
\be
C\int_0^\delta\!\!\!\int_{\Omega }  \lambda ( s\theta )^2  |v|^2dxdt \le -I  + \frac{1}{2}\int_{\Omega }  [(1- \zeta )^2   (s\theta)  |v|^2 ]_{\vert t=\delta} dx 
+C'  \int_0^\delta\!\!\! \int_{ X(\omega ,t,0) } \lambda  ( s\theta ) ^2  |v|^2dxdt. \label{F1} 
\ee

Finally, by  Cauchy-Schwarz inequality, we have for any $\kappa >0$ 
\be
\vert I \vert  \le (4 \kappa )^{-1} \int_0^\delta \!\!\! \int_{\Omega}  |e^{-s\varphi } q_t |^2dxdt + \kappa \int_0^\delta \!\!\! \int_{\Omega }  ( s \theta  )^2  |v|^2dxdt. \label{F2} 
\ee
Combining \eqref{F1} with \eqref{F2}  gives \eqref{F0} for $\kappa /\lambda  >0$ small enough. Therefore, Claim 3 is proved. \qed

We can prove in the same way the following estimate for $t\in [T-\delta, T]$:
\begin{multline}
\int_{T-\delta}^T \!\int_{\Omega } \lambda  ( s\theta ) ^2  | v |^2  dxdt 
\le 
C\left( \int_{T-\delta}^T \!\ \int_{\Omega } 
\lambda ^{-1} | e^{-s\varphi} q_t |^2 dxdt \right. \\
\left. + 
\int_{\Omega } [(1-\zeta ) ^2 ( s\theta ) |v|^2]_{\vert t= T- \delta }dx + \int_{T-\delta}^T   \!\int_{ X(\omega ,t,0) } \lambda  ( s\theta ) ^2 
 |v|^2 dxdt   \right). \label{F0bis} 
\end{multline}
Let us now consider the estimate for $t\in [\delta , T-\delta ]$. \\
{\sc Claim 4.}
\begin{multline}
\int_\delta^{T-\delta} \!\!\! \int_{\Omega }  \lambda ^2 ( s\theta ) |v|^2dxdt + \int_{\Omega }  [(1-\zeta )^2(\lambda  s\theta) |v|^2]_{\vert t=\delta} dx 
+ \int_{\Omega } [(1-\zeta )^2(\lambda  s\theta) |v|^2]_{\vert t=T-\delta } dx \\
\le C \left( \int_\delta ^{T-\delta } \!\!\! \int_{\Omega } |e^{-s\varphi} q_t|^2 dxdt + \int_\delta ^{T-\delta}  \!\!\! \int_{ X(\omega ,t,0) } \lambda ^2 ( s\theta ) |v|^2 dxdt \right) . 
\label{F1000}
\end{multline} 
$\Vert\cdot\Vert$ and $(.,.)$ denoting here the Euclidean norm and scalar product in $L^2(\Omega \times (\delta ,T - \delta ))$, we have that 
\be
\label{XYZ0}
||e^{-s\varphi} q_t ||^2 \ge || (1-\zeta ) (s\varphi _ t  v +  v_t ) ||^2 \ge 2 ( (1-\zeta ) s\varphi _t v, (1-\zeta ) v_t). 
\ee
Next, we compute 
\begin{multline}
((1-\zeta )s \varphi _t v,(1-\zeta ) v_t ) 
= \int_{\Omega } (1-\zeta  )^2 s \varphi _t  \frac{v^2}{2}dx \bigg\vert _{t= \delta}^{T-\delta}  
-\frac{s}{2}\int_\delta^{T-\delta} \!\!\! \int_{\Omega } (1-\zeta  )^2 \varphi _{tt} |v|^2 dxdt  \\
 - \int_\delta^{T-\delta} \!\!\! \int_{\Omega } (1-\zeta) \nabla\xi ( X(x,0,t) ) \cdot 
  \big( \frac{\partial X}{\partial x} \big) ^{-1} (X(x,0,t),t,0) f(x,t)  s  \varphi _t v^2 dxdt . \qquad \qquad  \label{XYZ1}
\end{multline}
Since $g(t)=1$ for $t\in [\delta , T-\delta ]$, we have that $\varphi _t = -\lambda \psi _t e^{\lambda \psi}$. From \eqref{P3}-\eqref{P4}, we infer that
\begin{eqnarray*}
s \varphi _t (x,T-\delta ) &\ge& C\lambda s e^{\lambda \psi} \qquad x\in \overline{\Omega} \setminus X(\omega _1, T-\delta , 0),\\
-s\varphi _t (x, \delta )   &\ge& C\lambda s e^{\lambda \psi} \qquad x\in \overline{\Omega} \setminus X(\omega _1,\delta ,0).  
\end{eqnarray*} 
Therefore, using \eqref{AW3},  
\be
\int_{\Omega } (1-\zeta )^2 s \varphi _t \frac{v^2}{2} dx \bigg\vert _\delta ^{T-\delta}  \ge 
C\left( 
\int_{\Omega }[(1-\zeta ) ^2  (\lambda s\theta ) |v|^2 ]_{\vert t=\delta} dx + \int_{\Omega } [(1-\zeta )^2  ( \lambda s\theta  ) |v|^2 ]_{\vert t=T - \delta }dx 
\right) . 
\label{XYZ2}
\ee
Next, with 
$\varphi_{tt} = -\{ (\lambda \psi _t ) ^2  +\lambda \psi _{tt}     \} e^{ \lambda \psi} $ and \eqref{P2}, we obtain for $\lambda \ge \lambda _1>\lambda _0$
\be
-\frac{s}{2} \int_\delta^{T-\delta}  \!\!\! \int_{\Omega }(1-\zeta  )^2 \varphi _{tt}  |v|^2dxdt  \ge  C  \int_\delta^{T-\delta} \!\!\! \int_{\Omega} (1-\zeta ) ^2  \lambda ^2 s\theta |v|^2 dxdt. 
\label{XYZ3}
\ee
Finally 
\begin{multline}
\label{XYZ4}
 \left\vert  \int_\delta^{T-\delta} \!\!\! \int_{\Omega } (1-\zeta)  \nabla\xi ( X(x,0,t )  ) \cdot 
  \big( \frac{\partial X}{\partial x} \big) ^{-1} (X(x,0,t),t,0) f(x,t)
   s  \varphi _t v^2 dxdt \right\vert \\
 \le C \int_{\delta}^{T-\delta}  \!\!\! \int_{ X(\omega ,t,0) } \lambda s\theta |v|^2.
\end{multline}
Claim 4 follows from \eqref{XYZ0}-\eqref{XYZ4}. 

We infer from \eqref{F0}, \eqref{F0bis} and \eqref{F1000}  that for some constants $\lambda _1\ge \lambda _0$, $s_1\ge s_0$ and $C_1>0$ we have 
for all $\lambda \ge \lambda _1$ and all $s\ge s_1$
\be
\int_0^T\!\!\! \int_{\Omega }   \lambda ^2  (s\theta  )|v|^2dxdt \le
 C_1\left( 
 \int_0^T\!\!\! \int_{\Omega }  |e^{-s\varphi} q_t|^2 dxdt+ \int_0^T \!\!\! \int_{ X(\omega  ,t,0) } \lambda ^2  ( s\theta )^2 |v|^2 dxdt   
 \right) . \label{F200} 
\ee
Replacing $v$ by $e^{-s\varphi}q$ in \eqref{F200} gives at once \eqref{E50}. The proof of Lemma \ref{lem2} is complete.\qed

\section*{Acknowledgements}
This work  started when the second author was visiting the Basque Center for Applied Mathematics at Bilbao in 2011. The second author thanks that institution
for its hospitality and its support. The second author has also been supported by Agence Nationale de la Recherche, Project CISIFS,
grant ANR-09-BLAN-0213-02.
The first and third authors have been partially supported by Grant MTM2011-29306-C02-00 of the MICINN, Spain, the ERC Advanced Grant FP7-246775 NUMERIWAVES, ESF Research Networking Programme OPTPDE. The first named author has also been supported  by the Grant PI2010-04 of the Basque Government.

\end{document}